\documentclass[11pt,a4paper]{amsart}
\usepackage[T1]{fontenc}
\usepackage[utf8]{inputenc}
\usepackage{lmodern,microtype,mathrsfs,amsmath,amssymb,amsthm,mathtools,booktabs,array}
\UseMicrotypeSet[expansion]{basicmath}
\usepackage[margin=27mm]{geometry}
\usepackage{needspace,pgfplots}
\pgfplotsset{compat=1.18}
\usepackage[hidelinks]{hyperref}
\hypersetup{pdftitle={Hardy transitions on product lattices},pdfauthor={Carlos Lizama}}
\numberwithin{equation}{section}
\newtheorem{theorem}{Theorem}[section]
\newtheorem{proposition}[theorem]{Proposition}
\newtheorem{lemma}[theorem]{Lemma}
\newtheorem{corollary}[theorem]{Corollary}
\theoremstyle{remark}
\newtheorem{remark}[theorem]{Remark}
\newcommand{\N}{\mathbb N}
\newcommand{\Z}{\mathbb Z}
\newcommand{\R}{\mathbb R}
\newcommand{\E}{\mathcal E}
\newcommand{\M}{\mathcal M}
\newcommand{\D}{\mathcal D}
\newcommand{\A}{\mathcal A}
\newcommand{\eps}{\varepsilon}
\newcommand{\ip}[2]{\langle #1,#2\rangle}
\DeclareMathOperator{\Var}{Var}
\DeclareMathOperator{\Rea}{Re}
\allowdisplaybreaks[1]
\title[Hardy transitions on product lattices]{Hardy transitions on product lattices}
\author{Carlos Lizama}
\address{Departamento de Matem\'atica y Ciencia de la Computaci\'on, Universidad de Santiago de Chile, Santiago, Chile}
\email{carlos.lizama@usach.cl}
\thanks{Supported by the Agencia Nacional de Investigaci\'on y Desarrollo (ANID), Chile, through FONDECYT Regular, grant No.~1260057.}
\date{}
\begin{document}
\begin{abstract}
Boundary conditions can change the dimensional scale of a sharp functional inequality. We determine this change for inverse-square Hardy inequalities on the product lattices $\mathbb N^k\times\mathbb Z^{N-k}$. Their sharp constants satisfy $C_{N,k}\sim N^2\Psi(k/N)$ uniformly for $1\le k\le N$, where $\Psi(\alpha)\sim2\alpha$ at zero. The law gives linear growth for a fixed number of positive coordinates, every intermediate power for suitable sublinear choices, and quadratic growth for a positive proportion. We develop a method that connects exact uncertainty inequalities with Hardy inequalities through logarithmic localization, and combines this transfer with a block decomposition to cover all proportions. The transfer holds on weighted graphs with an explicit geometric error and controls the variational problem over arbitrary functions. A mixture of two Mathieu ground energies determines $\Psi$. At the fully positive endpoint, the variational problem defining $\kappa=\Psi(1)$ has a unique nondegenerate minimum and a unique positive energy--moment minimizing profile. We also prove strict comparison with the continuum coefficient in every sector for $N\ge9$, the optimal uniform dimension.
\end{abstract}
\maketitle

\section{Introduction}

Sharp functional inequalities measure the cost of concentrating a function under geometric constraints. Their best constants also determine stability thresholds for differential and difference operators. Understanding how those constants change under discretization is therefore a common problem in analysis and spectral theory. Dimension and boundary conditions are basic geometric parameters, yet their combined effect can be very different on a lattice and in a continuous space.

The Hardy inequality gives an explicit continuous reference. Let $N\ge1$ be the dimension and let $0\le k\le N$ count the coordinates constrained to be positive. On
\[
 K_{N,k}=(0,\infty)^k\times\R^{N-k}\quad(k\ge1),
 \qquad K_{N,0}=\R^N\setminus\{0\},
\]
it compares the Dirichlet energy $\int|\nabla u|^2$ with the weighted mass $\int|u|^2/|x|^2$, where $|x|$ is the Euclidean norm and the functions vanish on the coordinate faces. The sharp coefficient is
\[
 A_{N,k}=\frac{(N+2k-2)^2}{4}.
\]
For $k=0$ this is the classical Hardy constant; the formula for positive $k$ is due to Su and Yang \cite{SuYang}. This single formula gives quadratic growth in dimension throughout the family.

On the corresponding lattices, write $\N=\{1,2,\ldots\}$ and set
\[
 \Omega_{N,k}=\N^k\times\Z^{N-k}\quad(k\ge1),
 \qquad\Omega_{N,0}=\Z^N\setminus\{0\}.
\]
Let $C_{N,k}$ be the sharp inverse-square Hardy coefficient for the nearest-neighbor energy with zero extension across the boundary. Gupta \cite{GuptaLattice} established linear growth on the punctured lattice. Its leading coefficient was subsequently identified by Gupta \cite[Theorem~1, $\ell=1$]{GuptaAsymptotics} and independently by Huang and Ye \cite[Theorem~1.1, $m=1$]{HuangYeAsymptotics}:
\[
 C_{N,0}\sim2N\qquad(N\to\infty).
\]
The dimensions at which the continuous coefficient loses sharpness are also known at the two endpoints. On the punctured lattice, equality holds for $3\le N\le8$ and strict inequality holds for $N\ge9$ \cite{AlpayThree,AlpayLow,LizamaLow,LizamaTransition}. On the positive integer lattice, equality holds in dimensions one and two and strict inequality holds from dimension three \cite{LizamaPositive}.

These results raise a question about the whole family: how does imposing more coordinate boundaries change the size of the sharp constant? Answering it requires estimates that remain sharp when the number of constrained coordinates varies with dimension. The main obstacle is the inverse-square weight, which couples all coordinates even though the underlying energy separates. Product trial functions suggest a leading coefficient; identifying the sharp constant requires a matching lower bound for arbitrary functions.

We determine the transition uniformly over the entire family. A fixed number of positive coordinates gives linear growth, sublinear powers give every exponent between one and two, and a positive proportion gives quadratic growth. A single coefficient curve governs these regimes. At the fully positive endpoint, its value is a constant $\kappa$ selected by a unique nondegenerate Mathieu minimum. Thus coordinate boundaries change the dimensional exponent of the discrete inequality, while the continuous family retains the quadratic scale.

To prove this law, we develop a method that connects a separable uncertainty problem to the original Hardy problem. Logarithmic localization converts the sharp energy--moment inequality into a Hardy inequality with an explicit error governed by the geometry of edges. A block decomposition extends the resulting precision to vanishing proportions of constrained coordinates. Together these steps control the infimum over all functions and establish the leading coefficient uniformly. The transfer is formulated for weighted graphs, separating the spectral computation of an uncertainty constant from the geometric estimate needed to recover a Hardy constant.

\subsection{The uniform scaling law}
Every function on $\Omega_{N,k}$ is extended by zero to $\Z^N$. The Dirichlet lattice Laplacian is
\begin{equation}\label{eq:laplacian}
 (L_{N,k}u)(x)=2Nu(x)-\sum_{\substack{y\in\Omega_{N,k}\\|x-y|=1}}u(y).
\end{equation}
Edges joining the domain to its boundary contribute to the energy. For the finitely supported functions $C_c(\Omega_{N,k})$, define
\begin{equation}\label{eq:hardydef}
 C_{N,k}=\inf_{0\ne u\in C_c(\Omega_{N,k})}
 \frac{\E[u]}{\D[u]},\qquad
 \E[u]=\ip{u}{L_{N,k}u},\quad
 \D[u]=\sum_{x\in\Omega_{N,k}}\frac{|u(x)|^2}{|x|^2}.
\end{equation}
The inner product is that of $\ell^2(\Omega_{N,k})$. The same definitions apply to $k=0$, with zero value at the removed vertex.

Two one-dimensional ground energies determine the limiting coefficient. Let $L_+$ and $L_0$ be the Laplacians on $\ell^2(\N)$ and $\ell^2(\Z)$, respectively, with diagonal $2$ and nearest-neighbor entries $-1$, and impose the zero boundary value at the origin for $L_+$. If $X$ denotes multiplication by the coordinate, put
\begin{equation}\label{eq:eta}
 \eta_+(\lambda)=\min\sigma(L_++\lambda X^2),\qquad
 \eta_0(\lambda)=\min\sigma(L_0+\lambda X^2),\qquad \lambda>0.
\end{equation}
These are the ground energies of discrete harmonic oscillators. For a proportion $0<\alpha\le1$, define
\begin{equation}\label{eq:psidef}
 h_\alpha(\lambda)=\alpha\eta_+(\lambda)+(1-\alpha)\eta_0(\lambda),\qquad
 \Psi(\alpha)=\inf_{\lambda>0}\frac{h_\alpha(\lambda)^2}{4\lambda},
\end{equation}
and put $\Psi(0)=0$. Both sectors use the same optimizing parameter $\lambda$.

\begin{theorem}[Uniform scaling law]\label{thm:main}
The function $\Psi$ is positive on $(0,1]$, and
\begin{equation}\label{eq:uniform}
 C_{N,k}\sim N^2\Psi(k/N)+2N\mathbf1_{\{k=0\}}
 \qquad\text{uniformly for }0\le k\le N\text{ as }N\to\infty.
\end{equation}
\end{theorem}

Uniformity means that the relative error tends to zero for every choice of $k=k_N$. Proposition~\ref{prop:uniformrelative} bounds that error by $20N^{-1/3}$. The indicator in \eqref{eq:uniform} accounts for the punctured lattice, whose leading term is $2N$; for $k\ge1$, the coefficient is governed by $\Psi$. Its expansion $\Psi(\alpha)\sim2\alpha$ at zero gives the regimes in Table~\ref{tab:scales}, where $\kappa=\Psi(1)$.

\begin{table}[ht]
\caption{A choice of the number of constraints selects the dimensional scale. The second column gives a quantity asymptotic to $C_{N,k}$ as $N\to\infty$. Here $m\ge1$ and $c>0$ are fixed, $0<\beta<1$, and $0<\alpha\le1$.}\label{tab:scales}
\centering
\begin{tabular}{@{}ll@{}}
\toprule
Choice of $k=k_N$ & Leading term\\
\midrule
$k=0$ (punctured lattice) & $2N$\\
$k=m\ge1$ fixed & $2mN$\\
$k\sim cN^\beta$, $0<\beta<1$ & $2cN^{1+\beta}$\\
$k/N\to\alpha>0$ & $\Psi(\alpha)N^2$\\
$k=N$ & $\kappa N^2$\\
\bottomrule
\end{tabular}
\end{table}

Every fixed positive proportion gives quadratic growth. Intermediate exponents arise when $k/N$ tends to zero at a prescribed rate. The table identifies both the exponent and its leading coefficient for each such choice.

\subsection{The transfer method}
Odd reflection gives a common spectral setting for the products. Extending a function on $\N$ with opposite sign across the origin identifies it with an odd function on $\Z$. Repeating this operation identifies $\Omega_{N,k}$ with a prescribed reflection sector of the full lattice. We call the reverse operation folding; each fold corresponds to one coordinate boundary.

The method starts from the reciprocal of the Hardy weight. With
\[
 \A[u]=\|u\|_2^2,\qquad \M[u]=\sum_x|x|^2|u(x)|^2,
\]
separation of $L_{N,k}+\lambda|x|^2$ and optimization in $\lambda$ give the exact inequality
\[
 \E[u]\M[u]\ge N^2\Psi(k/N)\A[u]^2.
\]
To transfer this information to the reciprocal moment $\D[u]$, we localize in $\log|x|$. On each localized piece, Kantorovich's inequality compares the two moments, while a discrete localization identity controls the energy cost. This yields
\[
 0\le N^2\Psi(k/N)-C_{N,k}\le4\pi\sqrt{10}\,N
 \qquad(1\le k\le N).
\]
Proposition~\ref{prop:transfer} establishes the underlying estimate on weighted graphs: the error is controlled by the sharp uncertainty constant and the variation of the logarithmic radius along edges. For families in which the latter is negligible relative to the former, the two sharp constants are asymptotically equal.

The block decomposition completes the method at small proportions. A direct positive comparison gives $C_{N,0}=2N+O(\sqrt N)$ and supplies sharp leading estimates for products with one positive coordinate. Splitting the coordinates into such blocks yields $C_{N,k}\sim2Nk$ whenever $1\le k=o(N)$. Combining this estimate with the localization bound proves Theorem~\ref{thm:main}. In particular, correlations between coordinates cannot improve the leading coefficient in any of these regimes.

\subsection{The coefficient curve and its endpoint}
The one-dimensional oscillators describe a particle hopping between neighboring lattice sites in a quadratic potential. Fourier transformation gives the differential expression
\[
 -\lambda\frac{d^2}{d\theta^2}+2-2\cos\theta.
\]
After rescaling, this is Mathieu's operator. For $q>0$, let $a_0(q)$ be the lowest eigenvalue of
\begin{equation}\label{eq:mathieuop}
 -\frac{d^2}{dz^2}+2q\cos(2z)
\end{equation}
with periodic boundary conditions on a period of length $\pi$. Let $b_2(q)$ be its lowest eigenvalue on $(0,\pi/2)$ with Dirichlet boundary conditions, equivalently the lowest eigenvalue in the odd $\pi$-periodic sector. Folding selects between these two sectors.

\Needspace{9\baselineskip}
\begin{theorem}[The limiting curve]\label{thm:curve}
For $0<\alpha\le1$,
\begin{equation}\label{eq:mathieupsi}
 \Psi(\alpha)=\min_{q>0}
 \frac{[2q+(1-\alpha)a_0(q)+\alpha b_2(q)]^2}{16q}.
\end{equation}
The function $\Psi$ is strictly increasing, $\sqrt\Psi$ is concave on $[0,1]$, and $\Psi(\alpha)\sim2\alpha$ as $\alpha\downarrow0$.
\end{theorem}

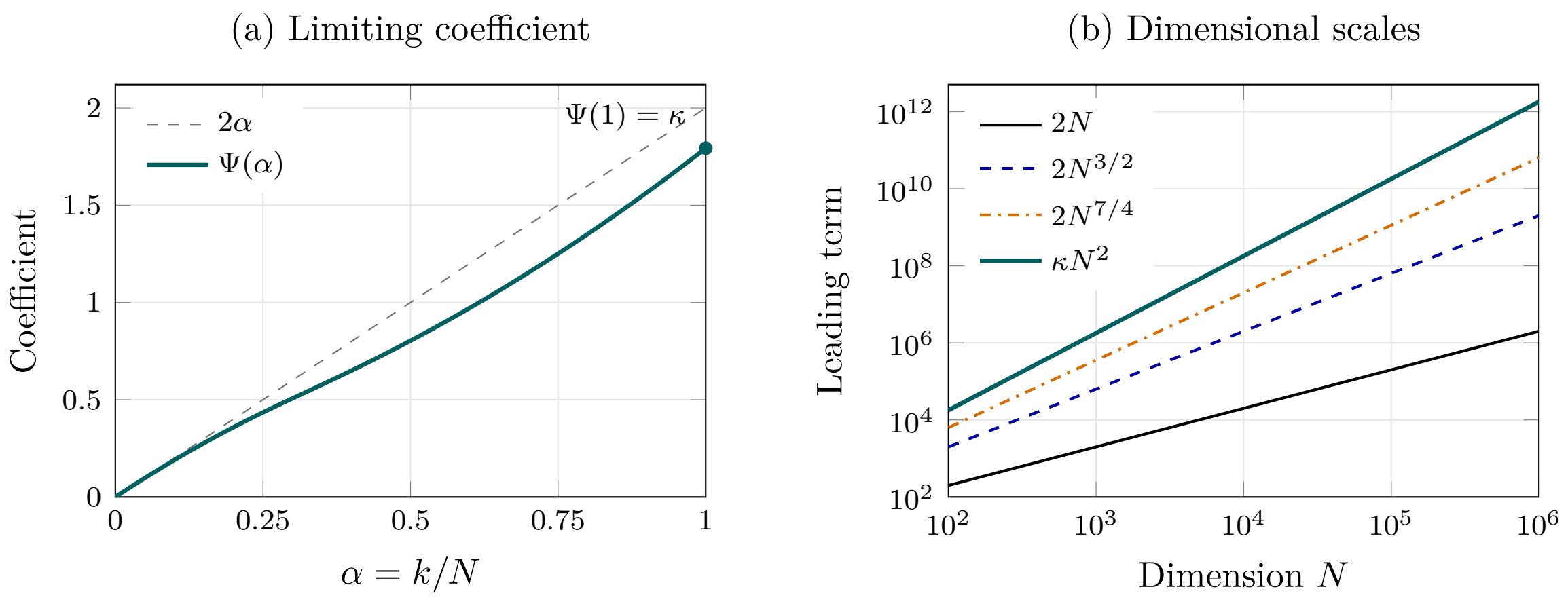
\begin{figure}[htbp]
\centering
\begin{tikzpicture}
\begin{axis}[
 name=profile,width=.47\textwidth,height=5.6cm,
 xmin=0,xmax=1,ymin=0,ymax=2.12,
 xlabel={$\alpha=k/N$},ylabel={Coefficient},
 title={(a) Limiting coefficient},
 xtick={0,.25,.5,.75,1},ytick={0,.5,1,1.5,2},
 tick label style={font=\scriptsize},label style={font=\small},
 title style={font=\small},grid=major,grid style={gray!18},
 legend style={font=\scriptsize,draw=none,at={(.03,.97)},anchor=north west},
 legend cell align=left]
\addplot[black!55,dashed,domain=0:1,samples=2] {2*x};
\addlegendentry{$2\alpha$}
\addplot[teal!75!black,very thick] coordinates {
(0.000,0.00000000000)
(0.005,0.00997497908)
(0.010,0.01989983187)
(0.015,0.02977443007)
(0.020,0.03959864314)
(0.025,0.04937233830)
(0.030,0.05909538051)
(0.035,0.06876763248)
(0.040,0.07838895467)
(0.045,0.08795920532)
(0.050,0.09747824044)
(0.055,0.10694591388)
(0.060,0.11636207739)
(0.065,0.12572658067)
(0.070,0.13503927145)
(0.075,0.14429999567)
(0.080,0.15350859754)
(0.085,0.16266491979)
(0.090,0.17176880386)
(0.095,0.18082009015)
(0.100,0.18981861836)
(0.105,0.19876422786)
(0.110,0.20765675812)
(0.115,0.21649604930)
(0.120,0.22528194282)
(0.125,0.23401428217)
(0.130,0.24269291378)
(0.135,0.25131768806)
(0.140,0.25988846070)
(0.145,0.26840509405)
(0.150,0.27686745898)
(0.155,0.28527543679)
(0.160,0.29362892177)
(0.165,0.30192782392)
(0.170,0.31017207242)
(0.175,0.31836161954)
(0.180,0.32649644532)
(0.185,0.33457656309)
(0.190,0.34260202598)
(0.195,0.35057293446)
(0.200,0.35848944533)
(0.205,0.36635178211)
(0.210,0.37416024717)
(0.215,0.38191523576)
(0.220,0.38961725208)
(0.225,0.39726692755)
(0.230,0.40486504127)
(0.235,0.41241254252)
(0.240,0.41991057492)
(0.245,0.42736050151)
(0.250,0.43476392957)
(0.255,0.44212273360)
(0.260,0.44943907415)
(0.265,0.45671540998)
(0.270,0.46395450064)
(0.275,0.47115939691)
(0.280,0.47833341723)
(0.285,0.48548010953)
(0.290,0.49260319987)
(0.295,0.49970653080)
(0.300,0.50679399416)
(0.305,0.51386946356)
(0.310,0.52093673178)
(0.315,0.52799945697)
(0.320,0.53506112040)
(0.325,0.54212499654)
(0.330,0.54919413496)
(0.335,0.55627135260)
(0.340,0.56335923441)
(0.345,0.57046014028)
(0.350,0.57757621633)
(0.355,0.58470940882)
(0.360,0.59186147954)
(0.365,0.59903402156)
(0.370,0.60622847471)
(0.375,0.61344614046)
(0.380,0.62068819573)
(0.385,0.62795570566)
(0.390,0.63524963525)
(0.395,0.64257085984)
(0.400,0.64992017455)
(0.405,0.65729830268)
(0.410,0.66470590323)
(0.415,0.67214357759)
(0.420,0.67961187545)
(0.425,0.68711130006)
(0.430,0.69464231295)
(0.435,0.70220533800)
(0.440,0.70980076521)
(0.445,0.71742895392)
(0.450,0.72509023577)
(0.455,0.73278491726)
(0.460,0.74051328210)
(0.465,0.74827559327)
(0.470,0.75607209488)
(0.475,0.76390301385)
(0.480,0.77176856135)
(0.485,0.77966893421)
(0.490,0.78760431608)
(0.495,0.79557487852)
(0.500,0.80358078201)
(0.505,0.81162217682)
(0.510,0.81969920383)
(0.515,0.82781199522)
(0.520,0.83596067519)
(0.525,0.84414536052)
(0.530,0.85236616115)
(0.535,0.86062318064)
(0.540,0.86891651666)
(0.545,0.87724626144)
(0.550,0.88561250205)
(0.555,0.89401532087)
(0.560,0.90245479583)
(0.565,0.91093100075)
(0.570,0.91944400559)
(0.575,0.92799387671)
(0.580,0.93658067710)
(0.585,0.94520446659)
(0.590,0.95386530207)
(0.595,0.96256323763)
(0.600,0.97129832478)
(0.605,0.98007061257)
(0.610,0.98888014773)
(0.615,0.99772697485)
(0.620,1.00661113647)
(0.625,1.01553267321)
(0.630,1.02449162385)
(0.635,1.03348802548)
(0.640,1.04252191357)
(0.645,1.05159332206)
(0.650,1.06070228342)
(0.655,1.06984882878)
(0.660,1.07903298796)
(0.665,1.08825478954)
(0.670,1.09751426095)
(0.675,1.10681142850)
(0.680,1.11614631746)
(0.685,1.12551895210)
(0.690,1.13492935573)
(0.695,1.14437755077)
(0.700,1.15386355877)
(0.705,1.16338740049)
(0.710,1.17294909587)
(0.715,1.18254866414)
(0.720,1.19218612381)
(0.725,1.20186149271)
(0.730,1.21157478805)
(0.735,1.22132602641)
(0.740,1.23111522377)
(0.745,1.24094239556)
(0.750,1.25080755668)
(0.755,1.26071072152)
(0.760,1.27065190396)
(0.765,1.28063111743)
(0.770,1.29064837489)
(0.775,1.30070368889)
(0.780,1.31079707155)
(0.785,1.32092853460)
(0.790,1.33109808940)
(0.795,1.34130574693)
(0.800,1.35155151784)
(0.805,1.36183541241)
(0.810,1.37215744063)
(0.815,1.38251761217)
(0.820,1.39291593640)
(0.825,1.40335242240)
(0.830,1.41382707899)
(0.835,1.42433991471)
(0.840,1.43489093785)
(0.845,1.44548015644)
(0.850,1.45610757831)
(0.855,1.46677321104)
(0.860,1.47747706197)
(0.865,1.48821913827)
(0.870,1.49899944688)
(0.875,1.50981799454)
(0.880,1.52067478783)
(0.885,1.53156983312)
(0.890,1.54250313660)
(0.895,1.55347470431)
(0.900,1.56448454212)
(0.905,1.57553265573)
(0.910,1.58661905070)
(0.915,1.59774373243)
(0.920,1.60890670619)
(0.925,1.62010797709)
(0.930,1.63134755014)
(0.935,1.64262543019)
(0.940,1.65394162197)
(0.945,1.66529613010)
(0.950,1.67668895907)
(0.955,1.68812011327)
(0.960,1.69958959697)
(0.965,1.71109741433)
(0.970,1.72264356942)
(0.975,1.73422806621)
(0.980,1.74585090854)
(0.985,1.75751210021)
(0.990,1.76921164487)
(0.995,1.78094954614)
(1.000,1.79272580751)
};
\addlegendentry{$\Psi(\alpha)$}
\addplot[only marks,mark=*,mark size=1.7pt,teal!75!black]
 coordinates {(1,1.79272580751)};
\node[font=\scriptsize,anchor=south east] at (axis cs:.99,1.81) {$\Psi(1)=\kappa$};
\end{axis}
\begin{axis}[
 at={(profile.right of south east)},anchor=left of south west,xshift=.75cm,
 width=.47\textwidth,height=5.6cm,
 xmin=2,xmax=6,ymin=2,ymax=12.7,
 xlabel={Dimension $N$},ylabel={Leading term},
 title={(b) Dimensional scales},
 xtick={2,3,4,5,6},xticklabels={$10^2$,$10^3$,$10^4$,$10^5$,$10^6$},
 ytick={2,4,6,8,10,12},yticklabels={$10^2$,$10^4$,$10^6$,$10^8$,$10^{10}$,$10^{12}$},
 tick label style={font=\scriptsize},label style={font=\small},
 title style={font=\small},grid=major,grid style={gray!18},
 legend style={font=\scriptsize,draw=none,at={(.03,.97)},anchor=north west},
 legend cell align=left,domain=2:6,samples=2]
\addplot[black,thick] {log10(2)+x};
\addlegendentry{$2N$}
\addplot[blue!65!black,thick,dashed] {log10(2)+1.5*x};
\addlegendentry{$2N^{3/2}$}
\addplot[orange!85!black,thick,dashdotted] {log10(2)+1.75*x};
\addlegendentry{$2N^{7/4}$}
\addplot[teal!75!black,very thick] {log10(1.79272580751)+2*x};
\addlegendentry{$\kappa N^2$}
\end{axis}
\end{tikzpicture}
\caption{Two aspects of the uniform law. (a) Numerical evaluation of \eqref{eq:mathieupsi}, together with the tangent $2\alpha$ at zero. (b) The leading terms for $k=1$, $k=\lfloor N^{1/2}\rfloor$, $k=\lfloor N^{3/4}\rfloor$, and $k=N$, respectively. Both axes in (b) are logarithmic, so the slopes are $1$, $3/2$, $7/4$, and $2$. The punctured case has the same leading term as $k=1$. The curves in (b) represent asymptotic leading terms. Computational details for (a) are given in Appendix~\ref{app:illustration}.}\label{fig:scaling}
\end{figure}

At $\alpha=1$, Theorem~\ref{thm:kappacharacterization} proves that the global minimum defining $\kappa$ is unique and nondegenerate and that its positive normalized energy--moment profile is unique. Sections~\ref{sec:smallalpha} and~\ref{sec:kappa} construct unique analytic minimizing branches near zero and one. In particular, the expansion at zero gives
\[
 C_{N,k}=2Nk-k^2+o(k^2)
 \quad\text{if }\sqrt N=o(k)\text{ and }k=o(N).
\]
Figure~\ref{fig:scaling} displays the coefficient curve and representative dimensional scales.

\subsection{Comparison with the continuum}
A dilation argument gives $C_{N,k}\le A_{N,k}$ for every pair $(N,k)$. The next theorem places every product in the strict range at a common dimension.

\begin{theorem}[Uniform strict comparison]\label{thm:comparison}
For every $N\ge9$ and every $0\le k\le N$,
\[
 C_{N,k}<A_{N,k}=\frac{(N+2k-2)^2}{4}.
\]
The dimension nine is optimal for this uniform assertion, since $C_{8,0}=A_{8,0}$.
\end{theorem}

Section~\ref{subsec:strictcomparison} proves strict comparison in the mixed sectors by explicit product trials and incorporates the punctured endpoint from \cite{LizamaTransition}. Table~\ref{tab:endpoints}, in that section, records the finite-dimensional picture and its remaining mixed cases. Corollary~\ref{cor:spectral} identifies $C_{N,k}$, up to the hopping coefficient, with the critical inverse-square coupling for the associated lattice Schr\"odinger operator.

\subsection{Relation to earlier work}
Continuous and discrete Hardy inequalities are compared in \cite{KapitanskiLaptev} and surveyed in \cite{Oberwolfach}; fractional and non-fractional lattice inequalities are studied in \cite{Dyda}. Ground-state representations \cite{FrankSeiringer} lead to optimal weights for elliptic operators \cite{DFP} and graphs \cite{KPP,KPPcriticality}. Keller and Lemm \cite{KellerLemm} determine the large-distance behavior of optimal lattice weights at fixed dimension. Here the inverse-square weight is prescribed and the dimension and coordinate constraints vary together.

The boundary theory includes weighted and higher-order half-line inequalities \cite{GuptaOneD,GKSRellich,HuangYeOneD} and optimal fractional weights \cite{DasFuente}. For the punctured lattice, higher-order leading asymptotics appear in \cite{GuptaAsymptotics,HuangYeAsymptotics} and finer corrections in \cite{LizamaAsymptotics}. Other symmetry constraints produce different dimensional laws: permutation antisymmetry gives fourth-order growth \cite{GuptaAntisymmetric}. Discrete rearrangement provides a further approach to lattice inequalities \cite{HajaiejHanHua}.

The spectral ingredients also occur in the discrete oscillator correspondence \cite{Aunola}, sharp number--angle uncertainty \cite{BKW}, and quantum estimation \cite{Hayashi}. The method developed here connects these one-dimensional spectral problems to sharp Hardy constants in varying dimension. Its quantitative transfer combines discrete localization \cite[\S3.1]{DiGesu} with Kantorovich's inequality \cite[Eq.~(1.1)]{ZhaoLiGong}; the block argument makes that connection uniform through the regime of vanishing boundary proportion.

Sections~\ref{sec:folding}--\ref{sec:uncertainty} establish the variational framework and exact uncertainty inequality. Sections~\ref{sec:localization}--\ref{sec:blocks} prove the uniform law, and Sections~\ref{sec:smallalpha}--\ref{sec:kappa} analyze the coefficient curve at its endpoints. Section~\ref{subsec:strictcomparison} proves the continuum comparison and its spectral consequence. The appendices provide the polynomial and perturbative calculations used in the proofs.

\section{Geometry and the variational framework}\label{sec:folding}

The boundary convention is part of the definition of the sharp constant. We first write the energy in a form that keeps every boundary contribution visible. We then explain folding and prove the continuum comparison. These preliminary facts also justify all truncations used later.

\subsection{Energy, domains, and finite truncations}
For this subsection, $\Omega$ is either $\Omega_{N,k}$ with $k\ge1$ or $\Z^N\setminus\{0\}$. Write $\widetilde u$ for the zero extension of $u$. Let $\mathscr E_{\rm int}$ be the set of unordered nearest-neighbor edges with both endpoints in $\Omega$, and let
\[
 b_\Omega(x)=\#\{y\notin\Omega:|y-x|=1\}.
\]
Thus $b_\Omega$ records the boundary contributions to the energy, with diagonal $2N$ in the Laplacian. All inner products below are conjugate-linear in the first variable.

\begin{lemma}\label{lem:forms}
The operator $L_{N,k}$ is bounded and self-adjoint on $\ell^2(\Omega)$, with $0\le L_{N,k}\le4NI$. Its energy is
\begin{align}
 \E[u]&=\sum_{i=1}^N\sum_{x\in\Z^N}
 |\widetilde u(x+e_i)-\widetilde u(x)|^2\label{eq:fullenergy}\\
 &=\sum_{\{x,y\}\in\mathscr E_{\rm int}}|u(x)-u(y)|^2
   +\sum_{x\in\Omega}b_\Omega(x)|u(x)|^2.\label{eq:edgeenergy}
\end{align}
Moreover, $\E[u]>0$ for $u\ne0$. Define
\[
 \mathscr Q_\Omega=\left\{u\in\ell^2(\Omega):
       \sum_{x\in\Omega}|x|^2|u(x)|^2<\infty\right\}.
\]
Finitely supported functions are dense in $\mathscr Q_\Omega$ for the norm
$(\|u\|_2^2+\M[u])^{1/2}$. Both $\E$ and $\D$ are continuous in $\ell^2(\Omega)$, and $\M$ is lower semicontinuous there. In particular, the infimum defining $C_{N,k}$ is unchanged if all nonzero functions in $\ell^2(\Omega)$ are allowed.
\end{lemma}
\begin{proof}
For finitely supported $u$, expanding the squares in \eqref{eq:fullenergy} gives $2N\sum|u|^2$ minus twice the real part of the sum over interior edges. This is $\ip{u}{L_{N,k}u}$. Separating interior edges from edges with one zero endpoint gives \eqref{eq:edgeenergy}. The adjacency operator has norm at most $2N$: this follows either from the row-sum bound or by compressing the sum of the $2N$ unitary shifts on the full lattice. Therefore $L_{N,k}$ is bounded and self-adjoint. Positivity follows from the energy formula, and
\[
 |\widetilde u(x+e_i)-\widetilde u(x)|^2
 \le2|\widetilde u(x+e_i)|^2+2|\widetilde u(x)|^2
\]
gives $\E[u]\le4N\|u\|_2^2$. Density extends the identities to every $u\in\ell^2$.

If $\E[u]=0$, every nearest-neighbor difference of $\widetilde u$ vanishes. Connectedness of $\Z^N$ makes $\widetilde u$ constant, and the only constant in $\ell^2(\Z^N)$ is zero. Thus $\E[u]>0$ for every nonzero $u\in\ell^2(\Omega)$.

For $u\in\mathscr Q_\Omega$, put $u_R(x)=u(x)\mathbf1_{\{|x|\le R\}}$. Both $\|u-u_R\|_2$ and $\M[u-u_R]$ tend to zero by convergence of the defining series. Boundedness of $L_{N,k}$ gives $\E[u_R]\to\E[u]$. Multiplication by $|x|^{-2}$ is bounded, since $|x|\ge1$ on every domain under consideration; for $k\ge1$ its norm is at most $1/k$. This proves continuity of $\D$. If $u_j\to u$ in $\ell^2$, then $u_j(x)\to u(x)$ at each vertex, so Fatou's lemma gives $\M[u]\le\liminf_j\M[u_j]$.

Finally, for any nonzero $u\in\ell^2$, the same finite truncations satisfy $\E[u_R]\to\E[u]$ and $\D[u_R]\to\D[u]>0$. Their Hardy quotients converge to that of $u$. The inequality between the two infima in the other direction is immediate from inclusion of the admissible classes.
\end{proof}

For general background on forms and operators on graphs, see \cite{KLW}.

\subsection{Folding by odd reflection}
Folding identifies a product domain with a symmetry subspace without changing any radial weighted quotient.

\begin{proposition}\label{prop:fold}
For $1\le k\le N$, define $U_k:\ell^2(\Omega_{N,k})\to\ell^2(\Z^N)$ by
\[
 (U_ku)(x)=2^{-k/2}\left(\prod_{j=1}^k\operatorname{sgn}(x_j)\right)
 u(|x_1|,\ldots,|x_k|,x_{k+1},\ldots,x_N)
\]
when $x_1\cdots x_k\ne0$, and by zero otherwise. Then $U_k$ is unitary onto the subspace of functions odd in each of the first $k$ coordinates. It intertwines $L_{N,k}$ with the full lattice Laplacian on that subspace and preserves radial weighted quadratic forms. For the weight $|x|^{-2}$, its value at the origin can be assigned arbitrarily in the full-space expression. Furthermore,
\begin{equation}\label{eq:monotonicity}
 C_{N,0}\le C_{N,1}\le\cdots\le C_{N,N},
 \qquad C_{N,k}\le2Nk\quad(1\le k\le N).
\end{equation}
\end{proposition}
\begin{proof}
Each $x\in\Omega_{N,k}$ has exactly $2^k$ reflected copies. The factor $2^{-k/2}$ therefore preserves the squared norm. Reflection preserves $|x|$, so the same calculation preserves every radial weighted form on its natural domain. Conversely, a function odd in the specified coordinates vanishes on their coordinate hyperplanes and is uniquely determined by its restriction to $\Omega_{N,k}$. Multiplying that restriction by $2^{k/2}$ gives the inverse map.

To check the operator identity, first consider a vertex away from all reflecting hyperplanes. Reflection sends its neighbors to the corresponding neighbors of the original vertex, with the prescribed signs, so the Laplacian commutes with $U_k$ there. At a vertex on exactly one reflecting hyperplane, the two normal neighbors have opposite values and cancel; all tangential neighbors have value zero. At an intersection of two or more reflecting hyperplanes, every neighbor remains on at least one of them and has value zero. Thus the full Laplacian also vanishes on the hyperplanes, and the intertwining identity holds at every vertex. Since the operators are bounded, checking it on finitely supported functions is sufficient.

For monotonicity, take a function supported on $\Omega_{N,k+1}$ and extend it by zero to $\Omega_{N,k}$. Its zero extension to $\Z^N$ is unchanged, so \eqref{eq:fullenergy} and the definition of $\D$ show that both numerator and denominator are unchanged. Taking infima proves $C_{N,k}\le C_{N,k+1}$. The embedding into $\Z^N\setminus\{0\}$ proves the first inequality in \eqref{eq:monotonicity}.

For the upper bound, let $x_*=(1,\ldots,1,0,\ldots,0)$ have $k$ unit coordinates. The function $\delta_{x_*}$ has energy $2N$ and weighted mass $|x_*|^{-2}=1/k$. Its quotient is $2Nk$.
\end{proof}

\subsection{The continuum coefficient}
The continuum coefficient of Su and Yang \cite{SuYang} provides a reference scale for every number of folds. We give the angular and radial arguments, including sharpness, to explain the discrete upper bound and fix the boundary convention.

\begin{proposition}\label{prop:continuum}
For $N\ge2$ and $0\le k\le N$, let
$K_{N,k}=(0,\infty)^k\times\R^{N-k}$ when $k\ge1$, and $K_{N,0}=\R^N\setminus\{0\}$. The sharp Hardy coefficient for functions in $C_c^\infty(K_{N,k})$ is
\begin{equation}\label{eq:continuum}
 A_{N,k}=\frac{(N+2k-2)^2}{4}.
\end{equation}
In addition, $C_{N,k}\le A_{N,k}$.
\end{proposition}
\begin{proof}
\emph{The angular eigenvalue.}
Let $\Sigma=K_{N,k}\cap\mathbb S^{N-1}$ and $P(x)=x_1\cdots x_k$, with the empty product $P=1$ when $k=0$.
Because each variable occurs only to its first power, $\Delta P=0$. Writing
$P(r\omega)=r^kp(\omega)$ in the polar-coordinate formula for the Laplacian gives
\[
 -\Delta_{\mathbb S^{N-1}}p=k(k+N-2)p.
\]
For $k\ge1$, the function $p$ is positive in $\Sigma$ and vanishes on $\partial\Sigma$. The spherical section is a Lipschitz domain, and the smooth function $p$ has zero trace, so $p\in H_0^1(\Sigma)$. For $k=0$, $\Sigma$ is the whole sphere and $p=1$ has eigenvalue zero; here $H_0^1(\Sigma)=H^1(\Sigma)$ because there is no angular boundary. The same calculation applies in both cases.
For $v\in C_c^\infty(\Sigma)$, integration by parts yields
\begin{equation}\label{eq:angularGS}
 \int_\Sigma|\nabla_{\mathbb S}v|^2
 -k(k+N-2)\int_\Sigma|v|^2
 =\int_\Sigma p^2\left|\nabla_{\mathbb S}(v/p)\right|^2\ge0.
\end{equation}
This identity proves the lower bound for the angular Rayleigh quotient.
To obtain the opposite bound, approximate $p$ in $H_0^1(\Sigma)$ by compactly supported smooth functions. Their Rayleigh quotients converge to $k(k+N-2)$.

\emph{The radial coefficient.}
Put $a=(N-2)/2$. For a smooth compactly supported radial factor $h$ on $(0,\infty)$, write $h(r)=r^{-a}w(\log r)$. The change of variables $t=\log r$ gives
\begin{align*}
 \int_0^\infty |h|^2r^{N-3}\,dr&=\int_\R|w|^2\,dt,\\
 \int_0^\infty |h'|^2r^{N-1}\,dr
 &=\int_\R|w'-aw|^2\,dt\\
 &=a^2\int_\R|w|^2\,dt+\int_\R|w'|^2\,dt.
\end{align*}
The mixed term vanishes because it is the integral of a derivative.
Applying this calculation for each angular variable, and applying
\eqref{eq:angularGS} for each radius, gives
\[
 \int_{K_{N,k}}|\nabla v|^2
 \ge\left(a^2+k(k+N-2)\right)
       \int_{K_{N,k}}\frac{|v|^2}{|x|^2}.
\]
The coefficient simplifies to \eqref{eq:continuum}.

For sharpness, choose angular functions $p_j\in C_c^\infty(\Sigma)$ whose Rayleigh quotients tend to $k(k+N-2)$. Fix a nonzero $\chi\in C_c^\infty(\R)$, set $w_R(t)=\chi(t/R)$, and take
$v_{j,R}(r,\omega)=r^{-a}w_R(\log r)p_j(\omega)$.
Its quotient is exactly
\[
 a^2+\frac{\int_\Sigma|\nabla_{\mathbb S}p_j|^2}{\int_\Sigma|p_j|^2}
 +\frac1{R^2}\frac{\int_\R|\chi'|^2}{\int_\R|\chi|^2}.
\]
Let $R\to\infty$ and then $j\to\infty$. Each test function has compact support inside the cone, so this proves sharpness in the asserted class.

\emph{Transfer to the lattice.}
Fix $v\in C_c^\infty(K_{N,k})$ and extend it smoothly by zero to $\R^N$.
Set $u_R(x)=v(x/R)$ on $\Omega_{N,k}$. For large $R$, all nonzero terms stay away from the zero boundary. Taylor expansion, uniformly on a fixed compact neighborhood of the support, gives
\[
 v((x+e_i)/R)-v(x/R)
 =R^{-1}\partial_i v(x/R)+O(R^{-2}).
\]
There are $O(R^N)$ relevant vertices. After summing the squared differences and multiplying by $R^{2-N}$, the error tends to zero. The corresponding Riemann sums give
\[
 R^{2-N}\E[u_R]\to\int_{K_{N,k}}|\nabla v|^2,\qquad
 R^{2-N}\D[u_R]\to\int_{K_{N,k}}\frac{|v|^2}{|x|^2}.
\]
The denominator converges to a positive number when $v\ne0$.
Thus $C_{N,k}$ does not exceed the continuum quotient of $v$. Taking its infimum completes the proof.
\end{proof}

For $N=1$, the continuum half-line coefficient is $1/4$, obtained by the same radial calculation with $a=-1/2$. The same dilation argument gives $C_{1,k}\le1/4$ for $k=0,1$.

The angular eigenvalue $k(k+N-2)$ and the radial Hardy coefficient combine into $A_{N,k}=(N+2k-2)^2/4$. Thus the continuum formula depends only on $N+2k$. When $k/N\to\alpha$, it gives the reference limit $A_{N,k}/N^2\to(1+2\alpha)^2/4$.

For finite-element approximation of the continuum Hardy constant, Ignat and Zuazua \cite{IgnatZuazua} determine a sharp logarithmic convergence rate as the mesh is refined. The comparison above instead uses dilated test functions on a fixed lattice.

\section{An exact uncertainty inequality}\label{sec:uncertainty}

Throughout this section, $1\le k\le N$. The reciprocal of the Hardy weight is a sum of coordinate squares. This makes the energy-moment problem separable, even though the original Hardy quotient is not. We establish the one-dimensional spectral facts first, then prove the exact uncertainty inequality for arbitrary functions, and finally identify its constant through Mathieu's equation.

\subsection{The one-dimensional oscillators}
For $I=\N$ or $\Z$, multiplication by the coordinate has domains
\[
 D(X)=\{f:\sum_{j\in I}j^2|f_j|^2<\infty\},\qquad
 D(X^2)=\{f:\sum_{j\in I}j^4|f_j|^2<\infty\}.
\]
On $\N$ we always set $f_0=0$ in the difference expression.
Let $L_I=2I-S-S^*$, where $S$ is the unilateral or bilateral shift as appropriate, and let
\[
 Q_{I,\lambda}=L_I+\lambda X^2,\qquad \lambda>0.
\]
The notation $D(X)$ will always refer to the relevant one-dimensional space.

\begin{lemma}\label{lem:oscillatorfacts}
The operator $Q_{I,\lambda}$ is self-adjoint on $D(X^2)$ and has form domain $D(X)$ and compact resolvent. Its lowest eigenvalue $\eta_I(\lambda)$ is positive and simple, with a unique positive normalized eigenvector $f_{I,\lambda}$. The eigenvalue and a locally chosen normalized real eigenvector depend real analytically on $\lambda>0$, and
\begin{equation}\label{eq:FH}
 \eta_I'(\lambda)=\|Xf_{I,\lambda}\|_2^2.
\end{equation}
For $I=\Z$, the ground eigenvector is even.
\end{lemma}
\begin{proof}
\emph{Self-adjointness and compactness.}
The diagonal operator $\lambda X^2$ is self-adjoint, and $L_I$ is a bounded self-adjoint perturbation with $0\le L_I\le4I$. The sum is therefore self-adjoint on $D(X^2)$. Its closed quadratic form is
\[
 q_{I,\lambda}[f]=\ip f{L_If}+\lambda\|Xf\|_2^2,\qquad f\in D(X).
\]
Since $L_I$ is bounded, the form norm is equivalent to
$(\|f\|_2^2+\|Xf\|_2^2)^{1/2}$ for fixed $\lambda>0$.
The embedding of this form domain into $\ell^2(I)$ is compact: on a bounded set in the form norm,
\[
 \sum_{|j|>R}|f_j|^2\le R^{-2}\|Xf\|_2^2
\]
is uniformly small, while the coordinates with $|j|\le R$ belong to a finite-dimensional space. Compactness of the form-domain embedding implies compactness of the resolvent.

Compact resolvent and nonnegativity imply that the spectral bottom is an eigenvalue. It is positive: an eigenvector at zero would make both nonnegative terms of the form vanish, contradicting strict positivity of the kinetic energy on nonzero $\ell^2$ vectors.

\emph{Positivity and simplicity.}
Replacing $f$ by $|f|$ leaves its norm and moment unchanged and cannot increase its kinetic energy. Thus a ground eigenvector can be taken nonnegative. If a nonnegative ground eigenvector vanishes at a vertex $j$, its eigenvalue recurrence gives
$f_{j-1}+f_{j+1}=0$, with the missing term omitted at the half-line boundary. Both neighbors must vanish. Repetition along the connected line or half-line makes every coordinate zero, a contradiction. The nonnegative ground vector is therefore strictly positive.

To prove simplicity, choose a real ground eigenvector. Its modulus is also a ground vector. The preceding argument shows that no coordinate can vanish. If two neighboring coordinates had opposite signs, taking the modulus would strictly lower the kinetic energy on that edge, which is impossible at the minimum. Hence every real ground vector has constant sign. If the real eigenspace had dimension at least two, it would contain a vector orthogonal to a positive ground vector; that vector could not have constant sign. The real eigenspace is one-dimensional, and the complex eigenspace is its complexification.

\emph{Parameter dependence and the derivative.}
Near any $\lambda_0>0$, the operators form an analytic family with common domain $D(X^2)$: for each vector in that domain, $Q_{I,\lambda}f$ depends affinely on $\lambda$. The simple isolated eigenvalue and its spectral projection are analytic by analytic perturbation theory \cite[Chapter VII]{Kato}. Fixing the sign of the positive eigenvector gives a locally analytic normalized choice. Differentiating the eigenvalue equation, taking its inner product with $f_{I,\lambda}$, and using self-adjointness cancels the eigenvector derivative and gives \eqref{eq:FH}. On $\Z$, reflection commutes with the operator; positivity and uniqueness imply that reflection fixes the ground vector.
\end{proof}

\subsection{Exact separation of energy and moment}
For $u\in\mathscr Q_{\Omega_{N,k}}$, write
\[
 \A[u]=\|u\|_2^2,\qquad \M[u]=\sum_x|x|^2|u(x)|^2.
\]
The elementary identity
\begin{equation}\label{eq:scalar}
 \frac{(E+\lambda M)^2}{4\lambda}
 =EM+\frac{(E-\lambda M)^2}{4\lambda}
 \quad(E,M,\lambda>0)
\end{equation}
shows that its infimum over $\lambda$ is $EM$, attained exactly at $\lambda=E/M$. This identity is the reason a product of two quadratic forms can be handled through a sum of operators.

\begin{proposition}\label{prop:uncertainty}
For $\rho=k/N$ and every $u\in\mathscr Q_{\Omega_{N,k}}$,
\begin{equation}\label{eq:uncertainty}
 \E[u]\M[u]\ge N^2\Psi(\rho)\A[u]^2.
\end{equation}
The coefficient is sharp, both on $\mathscr Q_{\Omega_{N,k}}$ and on finitely supported functions. For every real $\alpha\in(0,1]$, one also has
\begin{equation}\label{eq:profile}
 \Psi(\alpha)=\inf_{f,g}\bar e\,\bar m,
\end{equation}
where $f\in D(X)\subset\ell^2(\N)$ and $g\in D(X)\subset\ell^2(\Z)$ are normalized and
\begin{align*}
 e_+(f)&=\ip f{L_+f},&m_+(f)&=\|Xf\|_2^2,\\
 e_0(g)&=\ip g{L_0g},&m_0(g)&=\|Xg\|_2^2,\\
 \bar e&=\alpha e_+(f)+(1-\alpha)e_0(g),&
 \bar m&=\alpha m_+(f)+(1-\alpha)m_0(g).
\end{align*}
When $\alpha=1$, the factor $g$ is omitted.
\end{proposition}
\begin{proof}
\emph{The inequality for arbitrary functions.}
Begin with a finitely supported $u$ and fix a coordinate $i$. For each choice of the other coordinates, apply the one-dimensional oscillator inequality to the resulting section. The kinetic energy along coordinate $i$ includes its zero boundary term if $i\le k$. Summing over the remaining coordinates gives
\[
 \E_i[u]+\lambda\sum_xx_i^2|u(x)|^2
 \ge
 \begin{cases}
 \eta_+(\lambda)\A[u],&i\le k,\\
 \eta_0(\lambda)\A[u],&i>k.
 \end{cases}
\]
Here $\E_i$ is the $i$th summand in \eqref{eq:fullenergy}. Summing these $N$ inequalities yields
\begin{equation}\label{eq:sumoscillator}
 \E[u]+\lambda\M[u]\ge
 [k\eta_+(\lambda)+(N-k)\eta_0(\lambda)]\A[u]
 =Nh_\rho(\lambda)\A[u].
\end{equation}
For $u\ne0$, both $\E[u]$ and $\M[u]$ are positive. The two sides of \eqref{eq:sumoscillator} are nonnegative, so we may square and divide by $4\lambda$. Taking infima in $\lambda$ on both sides, and using \eqref{eq:scalar}, gives \eqref{eq:uncertainty}. Lemma~\ref{lem:forms} extends the inequality to every $u\in\mathscr Q_{\Omega_{N,k}}$.

\emph{The variational formula for the profiles.}
For any normalized pair $f,g$, the numbers $\bar e$ and $\bar m$ are positive; the latter is at least $\alpha$. Identity~\eqref{eq:scalar} gives
\[
 \bar e\,\bar m=\inf_{\lambda>0}
       \frac{(\bar e+\lambda\bar m)^2}{4\lambda}.
\]
The infimum over pairs and the infimum over $\lambda$ commute because both are infima over the same product set of triples $(f,g,\lambda)$. For fixed $\lambda$, the two profiles are independent, and Lemma~\ref{lem:oscillatorfacts} gives
\[
 \inf_{f,g}(\bar e+\lambda\bar m)
 =\alpha\eta_+(\lambda)+(1-\alpha)\eta_0(\lambda)=h_\alpha(\lambda).
\]
Since the quantities are nonnegative, taking the square commutes with this infimum. Therefore
\[
 \inf_{f,g}\bar e\,\bar m
 =\inf_{\lambda>0}\frac1{4\lambda}
       \left[\inf_{f,g}(\bar e+\lambda\bar m)\right]^2
 =\Psi(\alpha).
\]

\emph{Sharpness in dimension $N$.}
For $\alpha=\rho$, let
$u=f^{\otimes k}\otimes g^{\otimes(N-k)}$. Fubini's theorem and normalization give
\[
 \A[u]=1,\qquad \E[u]=k e_++(N-k)e_0=N\bar e,\qquad
 \M[u]=k m_++(N-k)m_0=N\bar m.
\]
Consequently, the infimum of $\E[u]\M[u]$ over such products is $N^2\Psi(\rho)$ by \eqref{eq:profile}. This gives the upper bound for the best constant in \eqref{eq:uncertainty}, matching the lower bound already proved.

To obtain the same sharpness using finite support, truncate $f$ and $g$ in their respective $D(X)$ norms and renormalize them. Their energies and moments converge, because the kinetic operators are bounded and the moments converge in the form norm. Their finite tensor products are finitely supported and have the same limiting energy-moment product.
\end{proof}

\begin{corollary}\label{cor:Hardyupper}
For every $1\le k\le N$,
\[
 C_{N,k}\le N^2\Psi(k/N).
\]
\end{corollary}
\begin{proof}
Cauchy--Schwarz gives $\A[u]^2\le\M[u]\D[u]$ for every nonzero finitely supported $u$. Thus $\E[u]/\D[u]\le\E[u]\M[u]/\A[u]^2$. Taking infima and using the sharpness assertion in Proposition~\ref{prop:uncertainty} proves the bound.
\end{proof}

The reciprocal-moment inequality $\D\M\ge\A^2$ gives the upper bound on the Hardy constant. For the lower bound, Section~\ref{sec:localization} obtains the reverse moment comparison on logarithmically localized pieces and controls the corresponding energy error.

\subsection{Fourier and sine transforms}
The transforms identify the two spectral sectors and fix every scaling factor. On $\ell^2(\Z)$ use
\[
 (\mathcal Ff)(\theta)=(2\pi)^{-1/2}\sum_{n\in\Z}f_n e^{in\theta},
 \qquad -\pi<\theta<\pi.
\]
The shift becomes multiplication by an exponential, and $X^2$ becomes $-d^2/d\theta^2$. Thus $\mathcal FQ_{\Z,\lambda}\mathcal F^{-1}$ is
\[
 -\lambda\frac{d^2}{d\theta^2}+2-2\cos\theta
\]
on the periodic $H^2$ domain. On $\ell^2(\N)$ use the unitary sine transform
\[
 (\mathcal Sf)(\theta)=\sqrt{\frac2\pi}\sum_{n\ge1}f_n\sin(n\theta),
 \qquad 0<\theta<\pi.
\]
The identity $2\cos\theta\sin(n\theta)=\sin((n+1)\theta)+\sin((n-1)\theta)$, with $\sin0=0$, gives the same potential. The transformed domain of $X^2$ is $H^2(0,\pi)\cap H_0^1(0,\pi)$.

For either operator, set $\theta=2z$ and factor out $\lambda/4$. The differential expression becomes
\[
 2+\frac{\lambda}{4}\left(-\frac{d^2}{dz^2}
                 -\frac8\lambda\cos(2z)\right).
\]
With $q=4/\lambda$, the potential in parentheses is initially $-2q\cos(2z)$. Translation by $\pi/2$ on the periodic interval changes its sign. On the Dirichlet interval $(0,\pi/2)$, reflection $z\mapsto\pi/2-z$ does the same and preserves both endpoints as a set. The lowest characteristic values are therefore
\begin{equation}\label{eq:mathetas}
 \eta_0(\lambda)=2+\frac{\lambda}{4}a_0(4/\lambda),\qquad
 \eta_+(\lambda)=2+\frac{\lambda}{4}b_2(4/\lambda).
\end{equation}
These unitary changes of variables identify both the operators and their boundary conditions. The discrete oscillator correspondence is discussed in \cite{Aunola}; the standard characteristic-value notation is that of \cite{DLMF}.
The same odd periodic sector occurs in Hayashi's quantum-estimation problem \cite[\S12.3]{Hayashi}, with spectral function
\[
 \gamma_{\mathrm{SU}(2)}(s)=1+\frac{s}{16}b_2(8/s)=\frac12\eta_+(s/2).
\]
Mathieu's equation also describes radio-frequency ion traps \cite{Paul,DLMF}; here its characteristic values determine the energy--moment balance of the discrete oscillator.

Substituting $q=4/\lambda$ into \eqref{eq:psidef} gives
\[
 \frac{h_\alpha(\lambda)^2}{4\lambda}
 =\frac{[2q+(1-\alpha)a_0(q)+\alpha b_2(q)]^2}{16q}.
\]
This proves the representation in Theorem~\ref{thm:curve}, initially with an infimum. The next lemma proves attainment.

\subsection{Attainment and uniform compactness of the parameters}
The two ends of the parameter interval have different meanings. At $\lambda\downarrow0$, the oscillators approach their continuum coefficients. At $\lambda\to\infty$, a positive proportion of half-line coordinates makes the quotient diverge. A finite-support trial profile lies strictly below the first endpoint value.

\begin{lemma}\label{lem:attainment}
For every $\alpha\in(0,1]$, the infimum in \eqref{eq:psidef} is attained, and
\begin{equation}\label{eq:strict}
 0<\Psi(\alpha)<2\alpha\le\frac{(1+2\alpha)^2}{4}.
\end{equation}
For each fixed $\alpha_0\in(0,1]$, there are a relative neighborhood $V$ of $\alpha_0$ in $(0,1]$ and a compact interval $K\subset(0,\infty)$ such that every minimizing parameter for every $\alpha\in V$ belongs to $K$.
\end{lemma}
\begin{proof}
\emph{Behavior at the two endpoints.}
The large-$q$ expansions in \cite[Eq.~28.8.1]{DLMF} give
\[
 a_0(q)=-2q+2\sqrt q+O(1),\qquad
 b_2(q)=-2q+6\sqrt q+O(1).
\]
Inserting $q=4/\lambda$ into \eqref{eq:mathetas} yields
\begin{equation}\label{eq:smalleta}
 \eta_0(\lambda)=\sqrt\lambda+O(\lambda),\qquad
 \eta_+(\lambda)=3\sqrt\lambda+O(\lambda)
 \quad(\lambda\downarrow0).
\end{equation}
Thus, with $G_\alpha(\lambda)=h_\alpha(\lambda)^2/(4\lambda)$ and
$c_\alpha=(1+2\alpha)^2/4$,
\[
 G_\alpha(\lambda)=c_\alpha+O(\sqrt\lambda).
\]
The remainder is uniform for $\alpha\in[0,1]$, since $h_\alpha$ is a convex combination of two fixed functions. At infinity, $m_+(f)\ge1$ and $L_+\ge0$ imply $\eta_+(\lambda)\ge\lambda$. Hence
\begin{equation}\label{eq:largeG}
 G_\alpha(\lambda)\ge\alpha^2\lambda/4.
\end{equation}

\emph{A value below the continuum endpoint.}
In \eqref{eq:profile}, take
$g=\delta_0$ and $f=(\delta_1+\eps\delta_2)/\sqrt{1+\eps^2}$.
Their individual quantities are
\[
 e_0=2,\quad m_0=0,\quad
 e_+=2-\frac{2\eps}{1+\eps^2},\quad
 m_+=\frac{1+4\eps^2}{1+\eps^2}.
\]
Consequently,
\begin{equation}\label{eq:twopoint}
 \bar e\,\bar m=
 \left(2-\frac{2\alpha\eps}{1+\eps^2}\right)
 \frac{\alpha(1+4\eps^2)}{1+\eps^2}
 =2\alpha-2\alpha^2\eps+O(\eps^2).
\end{equation}
For each fixed $\alpha>0$, a sufficiently small positive $\eps$ gives
$\Psi(\alpha)<2\alpha$. Since $c_\alpha-2\alpha=(\alpha-1/2)^2$, the infimum is strictly below $c_\alpha$, including when $\alpha=1/2$.

\emph{Attainment.}
Choose a value of $G_\alpha$ strictly below $c_\alpha$. The limit at zero and the divergence in \eqref{eq:largeG} exclude both ends from a sufficiently low sublevel set. Continuity then gives a minimum on a compact interval. Its value is positive by Lemma~\ref{lem:oscillatorfacts}.

\emph{Uniform localization near a fixed proportion.}
Fix $\alpha_0>0$ and choose $\lambda_0$ with
$G_{\alpha_0}(\lambda_0)<c_{\alpha_0}$. Take $\gamma>0$ so that the gap is greater than $4\gamma$. By continuity, after restricting to a neighborhood $V$ of $\alpha_0$ we have $\alpha\ge\alpha_0/2$ and
\[
 G_\alpha(\lambda_0)\le c_\alpha-3\gamma\qquad(\alpha\in V).
\]
The uniform expansion at zero gives an $\ell>0$ such that
$G_\alpha(\lambda)\ge c_\alpha-\gamma$ for $0<\lambda<\ell$ and all $\alpha\in V$. No minimum can occur there. The values $G_\alpha(\lambda_0)$ are bounded on a smaller closed neighborhood. The uniform bound
$G_\alpha(\lambda)\ge\alpha_0^2\lambda/16$ then excludes $\lambda>R$ for a sufficiently large $R$. Thus $K=[\ell,R]$ contains all the minima under consideration.
\end{proof}

The compactness assertion also proves continuity of $\Psi$ on $(0,1]$ directly: on a fixed compact parameter interval, $G_\alpha$ is uniformly continuous in $\alpha$, so its minimum is continuous. It will later prevent minimizing branches from escaping when $\alpha$ approaches one.

\begin{lemma}\label{lem:balance}
Let $\lambda_\alpha$ be any minimizing parameter and take the normalized ground states $f_\alpha,g_\alpha$ of the two oscillators at that same parameter. Then
\begin{equation}\label{eq:balance}
 2\lambda_\alpha h_\alpha'(\lambda_\alpha)=h_\alpha(\lambda_\alpha),
 \qquad \bar e=\lambda_\alpha\bar m,\qquad
 \bar e\,\bar m=\Psi(\alpha).
\end{equation}
For $\alpha=1$, only $f_\alpha$ is needed.
\end{lemma}
\begin{proof}
The minimum is interior and $h_\alpha$ is differentiable. Differentiating $G_\alpha$ gives
\[
 G_\alpha'(\lambda)=
 \frac{h_\alpha(\lambda)}{4\lambda^2}
       \bigl(2\lambda h_\alpha'(\lambda)-h_\alpha(\lambda)\bigr).
\]
Its first factor is positive. Hence stationarity gives the first identity.
By \eqref{eq:FH}, $h_\alpha'=\bar m$ for the chosen profiles, and their eigenvalue equations give $h_\alpha=\bar e+\lambda\bar m$.
The stationary equation reduces to $\bar e=\lambda\bar m$.
Finally $G_\alpha=(\bar e+\lambda\bar m)^2/(4\lambda)=\bar e\,\bar m$ at that parameter.
\end{proof}

\subsection{Monotonicity and concavity}
The strict increase of the curve comes from the difference between the two symmetry sectors. Odd extension identifies the half-line oscillator with the odd subspace of the full-line oscillator. Since the full-line ground eigenvalue is simple and its eigenvector is even and strictly positive, compactness of the resolvent implies
\begin{equation}\label{eq:stricteta}
 \eta_+(\lambda)>\eta_0(\lambda)\qquad(\lambda>0).
\end{equation}

\begin{proposition}\label{prop:shape}
The extension $\Psi(0)=0$ is continuous on $[0,1]$. The function $\Psi$ is strictly increasing, $\sqrt\Psi$ is concave, and $\Psi$ is locally Lipschitz on $(0,1)$. With $\kappa=\Psi(1)$, one has
\begin{equation}\label{eq:curvebounds}
 \kappa\alpha^2\le\Psi(\alpha)<2\alpha
 \le\frac{(1+2\alpha)^2}{4}\qquad(0<\alpha\le1).
\end{equation}
\end{proposition}
\begin{proof}
At $\alpha=0$, the same spectral infimum equals zero: the trial vector $\delta_0$ gives $\eta_0(\lambda)\le2$, and hence
$0\le\inf_{\lambda>0}\eta_0(\lambda)^2/(4\lambda)\le\inf_{\lambda>0}1/\lambda=0$.
Therefore, on the entire interval $[0,1]$,
\begin{equation}\label{eq:concave}
 R(\alpha):=\sqrt{\Psi(\alpha)}
 =\inf_{\lambda>0}
 \frac{\eta_0(\lambda)+\alpha[\eta_+(\lambda)-\eta_0(\lambda)]}{2\sqrt\lambda}.
\end{equation}
As an infimum of affine functions of $\alpha$, the function $R$ is concave.

For strict monotonicity, let $\alpha<\beta$ and choose a minimizing parameter $\lambda_\beta$, whose existence is known because $\beta>0$. Equation~\eqref{eq:stricteta} gives
\[
 R(\alpha)\le\frac{h_\alpha(\lambda_\beta)}{2\sqrt{\lambda_\beta}}
 <\frac{h_\beta(\lambda_\beta)}{2\sqrt{\lambda_\beta}}
 =R(\beta).
\]

A finite concave function is continuous and locally Lipschitz in the interior of its interval. For example, its secant slopes on a compact interior interval are bounded by secant slopes joining that interval to two fixed exterior points. This applies to $R$, and hence to $R^2=\Psi$, since $R$ is bounded. At zero, $\Psi(\alpha)\le2\alpha$ gives continuity. At one, concavity gives $R(\alpha)\ge\alpha R(1)$, whereas monotonicity gives $R(\alpha)\le R(1)$. These bounds prove continuity there and yield $\Psi(\alpha)\ge\kappa\alpha^2$. The remaining bounds are those in Lemma~\ref{lem:attainment}.
\end{proof}

\section{From uncertainty to Hardy by logarithmic localization}\label{sec:localization}

The passage from uncertainty to Hardy depends on the variation of the radius along edges. We first formulate this step for a weighted graph, combining discrete localization \cite[\S3.1]{DiGesu} with the classical Kantorovich inequality. The lattice estimate then follows by bounding that variation and inserting the exact uncertainty constant.

\subsection{A transfer estimate on weighted graphs}
Let $V$ be a nonempty countable, locally finite graph with symmetric edge weights $w_{xy}>0$. Each unordered edge is counted once. For a function $b:V\to[0,\infty)$, set
\[
 \E[u]=\sum_{\{x,y\}}w_{xy}|u(x)-u(y)|^2+\sum_x b(x)|u(x)|^2.
\]
The second term includes possible Dirichlet boundary contributions. Given $r:V\to(0,\infty)$, define
\[
 \A[u]=\sum_x|u(x)|^2,\qquad
 \M[u]=\sum_x r(x)^2|u(x)|^2,\qquad
 \D[u]=\sum_x r(x)^{-2}|u(x)|^2.
\]
All infima below are over nonzero finitely supported functions on $V$.

\begin{proposition}[Uncertainty-to-Hardy transfer]\label{prop:transfer}
Put
\[
 K=\inf_u\frac{\E[u]\M[u]}{\A[u]^2},\qquad
 C=\inf_u\frac{\E[u]}{\D[u]},\qquad
 \Gamma=\sup_{x\in V}r(x)^2\sum_{y\sim x}w_{xy}
                 (\log r(x)-\log r(y))^2.
\]
If $K<\infty$ and $\Gamma<\infty$, then, for every $\delta>0$,
\begin{equation}\label{eq:transfer}
 \frac{K}{\cosh^2(2\delta)}-\frac{\pi^2\Gamma}{8\delta^2}
 \le C\le K,\qquad
 0\le K-C\le\pi\sqrt{2\Gamma K}.
\end{equation}
\end{proposition}
\begin{proof}
\emph{Step 1: localize in the logarithm of the radius.}
Fix a finitely supported $u$. For $\delta>0$, let
\[
 \varphi(s)=\begin{cases}\cos(\pi s/2),&|s|\le1,\\0,&|s|>1,\end{cases}
 \qquad \chi_t(x)=\delta^{-1/2}
         \varphi\!\left(\frac{\log r(x)-t}{\delta}\right),
 \qquad v_t=\chi_tu.
\]
Since $\|\varphi\|_2^2=1$ and $\|\varphi'\|_2^2=\pi^2/4$,
\begin{equation}\label{eq:partitionD}
 \int_\R\chi_t(x)^2\,dt=1,\qquad
 \int_\R\A[v_t]\,dt=\A[u],\qquad
 \int_\R\D[v_t]\,dt=\D[u].
\end{equation}
Expanding the square on each edge gives the integrated localization identity
\begin{equation}\label{eq:ims}
 \int_\R\E[v_t]\,dt-\E[u]
 =\sum_{\{x,y\}}w_{xy}\Rea(u(x)\overline{u(y)})
      \int_\R(\chi_t(x)-\chi_t(y))^2\,dt.
\end{equation}
The diagonal terms cancel by \eqref{eq:partitionD}, including those containing $b(x)$; the remaining cross terms give \eqref{eq:ims}. Local finiteness and finite support make all edge sums finite.

The translation estimate
$\|\varphi(\cdot+h)-\varphi\|_2\le |h|\|\varphi'\|_2$
follows by integrating $\varphi'$ and applying Cauchy--Schwarz. After scaling it gives
\[
 \int_\R(\chi_t(x)-\chi_t(y))^2\,dt
 \le\frac{\pi^2}{4\delta^2}(\log r(x)-\log r(y))^2.
\]
Using $\Rea(u(x)\overline{u(y)})\le(|u(x)|^2+|u(y)|^2)/2$
in \eqref{eq:ims} and regrouping by vertices yields
\begin{equation}\label{eq:imserror}
 \begin{split}
 \int_\R\E[v_t]\,dt-\E[u]
 &\le\frac{\pi^2}{8\delta^2}\sum_x|u(x)|^2
       \sum_{y\sim x}w_{xy}(\log r(x)-\log r(y))^2\\
 &\le\frac{\pi^2\Gamma}{8\delta^2}\D[u].
 \end{split}
\end{equation}

\emph{Step 2: compare the reciprocal moments by Kantorovich.}
For $v_t\ne0$, its support satisfies $e^{t-\delta}\le r(x)\le e^{t+\delta}$.
Apply the classical Kantorovich inequality \cite[Eq.~(1.1)]{ZhaoLiGong}
to the positive diagonal matrix with entries $r(x)^2$ on that finite support.
Its spectral bounds $a=e^{2(t-\delta)}$ and $d=e^{2(t+\delta)}$ give
\begin{equation}\label{eq:annularmoment}
 \M[v_t]\D[v_t]\le\frac{(a+d)^2}{4ad}\A[v_t]^2
                         =\cosh^2(2\delta)\A[v_t]^2.
\end{equation}
Thus the definition of $K$ implies
\[
 \E[v_t]\ge K\frac{\A[v_t]^2}{\M[v_t]}
       \ge\frac{K}{\cosh^2(2\delta)}\D[v_t]
\]
whenever $v_t\ne0$; the final inequality also holds for $v_t=0$.
Integrating and using \eqref{eq:partitionD}--\eqref{eq:imserror}, then taking the infimum of $\E[u]/\D[u]$, proves the lower bound in \eqref{eq:transfer}.

\emph{Step 3: optimize the width.}
Cauchy--Schwarz gives $\M[u]\D[u]\ge\A[u]^2$, hence
$\E[u]/\D[u]\le\E[u]\M[u]/\A[u]^2$ and $C\le K$.
Writing $B=\pi^2\Gamma/8$ and using $\tanh(2\delta)\le2\delta$, we obtain
\[
 0\le K-C\le K\tanh^2(2\delta)+B\delta^{-2}
                   \le4K\delta^2+B\delta^{-2}.
\]
For $K,B>0$, take $\delta^2=\sqrt B/(2\sqrt K)$; the right-hand side becomes
$4\sqrt{KB}=\pi\sqrt{2\Gamma K}$.
If $K=0$, nonnegativity and $C\le K$ give $C=0$; if $\Gamma=0$, let $\delta\downarrow0$ in the lower bound. This proves all cases.
\end{proof}

For a family with $K>0$, the condition $\Gamma/K\to0$ therefore implies $C/K\to1$. The uncertainty constant and the variation of $\log r$ can be estimated separately; product structure is used below to compute the former.

\subsection{Application to the product lattices}
Return to $\Omega_{N,k}$ with $1\le k\le N$, unit edge weights, $r(x)=|x|$, and $b=b_\Omega$ from \eqref{eq:edgeenergy}. The boundary enters through $|x|^2\ge k$.

\begin{lemma}\label{lem:geometry}
For every $x\in\Omega_{N,k}$,
\begin{equation}\label{eq:geometry}
 \sum_{\substack{y\in\Omega_{N,k}\\|x-y|=1}}
 (\log|x|-\log|y|)^2\le\frac{32+8N/k}{|x|^2}.
\end{equation}
\end{lemma}
\begin{proof}
Fix an interior neighbor $y$ of $x$, and put $s=|x|^2$ and $s'=|y|^2$. Both radii are at least one and differ by at most one, so
$|x|\le|y|+1\le2|y|$ and $\min(s,s')\ge s/4$. The mean value estimate gives
\[
 |\log|x|-\log|y||=\tfrac12|\log s-\log s'|
 \le\frac{|s-s'|}{2\min(s,s')}\le\frac{2|s-s'|}{s}.
\]
For the formal neighbors $x+e_i$ and $x-e_i$, the squared radii differ by $2x_i+1$ and $-2x_i+1$. Adding the nonnegative squared differences for neighbors outside the domain yields
\[
 \sum_{y\sim x}(\log|x|-\log|y|)^2
 \le\frac4{s^2}\sum_{i=1}^N[(2x_i+1)^2+(-2x_i+1)^2]
 =\frac{32}{s}+\frac{8N}{s^2}
 \le\frac{32+8N/k}{s}.
\]
Only interior radii occur inside logarithms.
\end{proof}

\begin{proposition}\label{prop:localization}
Let $\rho=k/N$, where $1\le k\le N$. For every $\delta>0$,
\begin{equation}\label{eq:finite}
 \frac{N^2\Psi(\rho)}{\cosh^2(2\delta)}
 -\frac{\pi^2(4+\rho^{-1})}{\delta^2}
 \le C_{N,k}\le N^2\Psi(\rho).
\end{equation}
Consequently, uniformly over $1\le k\le N$,
\begin{equation}\label{eq:rate}
 0\le N^2\Psi(\rho)-C_{N,k}
 \le4\pi N\sqrt{(4+\rho^{-1})\Psi(\rho)}
 \le4\pi\sqrt{10}\,N.
\end{equation}
In particular, if $k_N/N\to\alpha\in(0,1]$, then $C_{N,k_N}/N^2\to\Psi(\alpha)$.
\end{proposition}
\begin{proof}
Proposition~\ref{prop:uncertainty} identifies the sharp uncertainty constant as $K=N^2\Psi(\rho)$, including for finitely supported functions. Lemma~\ref{lem:geometry} gives $\Gamma\le32+8/\rho$. Inserting these values into Proposition~\ref{prop:transfer} proves \eqref{eq:finite} and the first error bound in \eqref{eq:rate}. Since $\Psi(\rho)\le2\rho$, one has $(4+\rho^{-1})\Psi(\rho)\le8\rho+2\le10$. Division by $N^2$ and continuity of $\Psi$ give the asserted limit.
\end{proof}

\begin{remark}[Products realizing the limit]
Fix $\alpha\in(0,1)$ and take the normalized ground profiles $f_\alpha,g_\alpha$ in Lemma~\ref{lem:balance}, with averaged energy $\bar e$ and moment $\bar m$. If $k_N/N\to\alpha$, the products
\[
 u_N=f_\alpha^{\otimes k_N}\otimes g_\alpha^{\otimes(N-k_N)}
\]
satisfy $\E[u_N]/N\to\bar e$. Under $|u_N|^2$, the coordinate squares are independent with two fixed distributions of finite mean. The weak law of large numbers gives $|x|^2/N\to\bar m$ in probability. Since $N/|x|^2\le N/k_N$ is uniformly bounded, it follows that $N\D[u_N]\to\bar m^{-1}$. Thus their Hardy quotients divided by $N^2$ tend to $\bar e\bar m=\Psi(\alpha)$.

For $\alpha=1$, use the half-line minimizing profile and any fixed normalized line profile of finite second moment in the $o(N)$ remaining coordinates. The latter contribute $o(1)$ to $|x|^2/N$ in probability, by Markov's inequality, and the same conclusion follows. Lemma~\ref{lem:forms} supplies finitely supported approximations with these limits.
\end{remark}

\section{The puncture and a vanishing proportion of folds}\label{sec:blocks}

The absolute error in \eqref{eq:rate} has the same order as the leading constant when $k$ is fixed. The punctured case also lies outside the separated oscillator mixture: at $\alpha=0$, that mixture allows concentration at the origin and has value zero. We first account for the deleted vertex by a positive comparison function. Its estimate supplies the one-fold inequalities needed for a block decomposition and gives a quantitative relative law for every sector.

\subsection{A positive comparison on the punctured lattice}
The limit $C_{N,0}/N\to2$ was established by Gupta \cite[Theorem~1, $\ell=1$]{GuptaAsymptotics} and independently by Huang and Ye \cite[Theorem~1.1, $m=1$]{HuangYeAsymptotics}. We give a direct proof for the first-order constant, with an explicit remainder. It uses the following elementary form of the ground-state identity. The comparison function need not belong to $\ell^2$.

\begin{lemma}\label{lem:positivecomparison}
Let $h:\Omega_{N,0}\to(0,\infty)$, extend it by $h(0)=0$, and interpret $L_{N,0}h$ pointwise. For every finitely supported $u$,
\begin{equation}\label{eq:discreteGS}
 \E[u]-\sum_{x\in\Omega_{N,0}}
       \frac{(L_{N,0}h)(x)}{h(x)}|u(x)|^2
 =\sum_{\{x,y\}\text{ interior}}h(x)h(y)
       \left|\frac{u(x)}{h(x)}-\frac{u(y)}{h(y)}\right|^2.
\end{equation}
Consequently, if $|x|^2(L_{N,0}h)(x)/h(x)\ge c$ at every vertex, then $C_{N,0}\ge c$.
\end{lemma}
\begin{proof}
Only edges meeting the finite support of $u$ contribute. Expanding the square on the right gives
\[
 \sum_x\left(\sum_{\substack{y\in\Omega_{N,0}\\y\sim x}}
                 \frac{h(y)}{h(x)}\right)|u(x)|^2
       -2\sum_{\{x,y\}\text{ interior}}\Rea(\overline{u(x)}u(y)).
\]
The coefficient in parentheses is $2N-(L_{N,0}h)(x)/h(x)$. Substitution in the energy formula proves the identity. In particular, the diagonal $2N$ retains the contributions of the edges incident to the removed vertex. Nonnegativity of the right side and the pointwise lower bound give $\E[u]\ge c\D[u]$. Taking infima proves the last assertion.
\end{proof}

\begin{proposition}[The punctured-lattice estimate]\label{prop:punctured}
For every integer $N\ge1$,
\begin{equation}\label{eq:puncturedbound}
 2N-\max\{2\sqrt{6N},32\}\le C_{N,0}\le2N.
\end{equation}
In particular,
\begin{equation}\label{eq:known}
 C_{N,0}=2N+O(\sqrt N),\qquad C_{N,0}\sim2N.
\end{equation}
\end{proposition}
\begin{proof}
\emph{Step 1: choose the comparison function.}
Write $s=|x|^2$ and let $a>0$. Define
\[
 h_a(x)=
 \begin{cases}
 a,&s=1,\\
 (s-1)^{-1},&s\ge2,
 \end{cases}
 \qquad h_a(0)=0.
\]
It is positive on the punctured lattice. We estimate
$W_a(x)=s(L_{N,0}h_a)(x)/h_a(x)$ uniformly in $x$. The separate value on the first sphere removes the pole of $(s-1)^{-1}$ and will balance the costs of the adjacent layers.

\emph{Step 2: estimate every vertex with $s\ge5$.}
For these vertices, every neighbor has squared radius at least two. Indeed, an integer neighbor of a unit vector has squared radius zero, two, or four. Thus the rational branch of $h_a$ applies to every neighbor. Since $x_i^2\le s$ and $s>4$, all denominators below are positive, and
\[
 \frac{h_a(x+e_i)+h_a(x-e_i)}{h_a(x)}
 =(s-1)\left(\frac1{s+2x_i}+\frac1{s-2x_i}\right)
 =\frac{2s(s-1)}{s^2-4x_i^2}.
\]
Consequently,
\begin{align*}
 W_a(x)
 &=2Ns-2s^2(s-1)\sum_{i=1}^N\frac1{s^2-4x_i^2}\\
 &=2N-8(s-1)\sum_{i=1}^N\frac{x_i^2}{s^2-4x_i^2}.
\end{align*}
Using $s^2-4x_i^2\ge s(s-4)$ and $\sum_i x_i^2=s$, we obtain
\[
 W_a(x)\ge2N-\frac{8(s-1)}{s-4}\ge2N-32.
\]
The last fraction decreases for $s\ge5$ and equals $32$ at $s=5$. This estimate covers all coordinate configurations at once.

\emph{Step 3: check the remaining layers.}
The squared coordinates are integers in $\{0,1,4,9,\ldots\}$. Up to signed permutations, the only possibilities with $1\le s\le4$ are those in the following table. A row requiring more than $N$ nonzero coordinates is simply absent.
\begin{center}
\begin{tabular}{@{}lll@{}}
\toprule
Squared radius & Representative & $W_a(x)$\\\midrule
$1$ & $e_1$ & $2N-(2N-5/3)/a$\\
$2$ & $e_1+e_2$ & $2N+3-4a$\\
$3$ & $e_1+e_2+e_3$ & $2N-48/5$\\
$4$ & $2e_1$ & $2N+9/2-12a$\\
$4$ & $e_1+e_2+e_3+e_4$ & $2N-8$\\\bottomrule
\end{tabular}
\end{center}
Here are the neighbor counts determining the entries involving $a$. At $e_1$, the nonzero neighboring values sum to $1/3+2(N-1)$, whereas $h_a(e_1)=a$. At $e_1+e_2$, they sum to
\[
 2a+2\cdot\frac14+2(N-2)\cdot\frac12=2a+N-\frac32,
\]
and $h_a=1$. At $2e_1$, they sum to
$a+1/8+2(N-1)/4$, and $h_a=1/3$.
Multiplication of the corresponding Laplacian ratios by $s$ gives the three entries. For $s=3$, and for the four-unit-coordinate configuration with $s=4$, no neighbor lies on the first sphere. The rational identity from Step~2 applies, giving respectively
\[
 8(3-1)\frac3{3^2-4}=\frac{48}{5},\qquad
 8(4-1)\frac4{4^2-4}=8
\]
as the subtracted terms.

\emph{Step 4: balance the exceptional layers.}
Take $a=\sqrt{N/6}$. Then
\[
 \frac{2N-5/3}{a}\le\frac{2N}{a}=2\sqrt{6N},
 \qquad 4a-3\le12a=2\sqrt{6N},
 \qquad12a-\frac92\le2\sqrt{6N}.
\]
The other two corrections in the table are less than $32$. Combining all layers gives
\[
 |x|^2\frac{(L_{N,0}h_a)(x)}{h_a(x)}
 \ge2N-\max\{2\sqrt{6N},32\}
 \quad(x\ne0).
\]
Lemma~\ref{lem:positivecomparison} proves the lower bound. For the upper bound use $u=\delta_{e_1}$, whose energy is $2N$ and Hardy mass is one. Dividing the resulting bounds by $2N$ proves \eqref{eq:known}.
\end{proof}

Monotonicity and the point-mass trial now give the one-fold estimates needed for the block argument:
\begin{equation}\label{eq:onefold}
 2n-\max\{2\sqrt{6n},32\}\le C_{n,0}\le C_{n,1}\le2n.
\end{equation}
In particular $C_{n,1}\sim2n$, and $C_{n,1}\ge n$ for every $n\ge32$. 

\subsection{Combining one-fold blocks}
\begin{proposition}\label{prop:blocks}
Write $N=kn+r$, where $1\le k\le N$, $n\ge1$, and $0\le r<k$. Then
\begin{equation}\label{eq:blocks}
 \bigl[(k-r)\sqrt{C_{n,1}}+r\sqrt{C_{n+1,1}}\bigr]^2
 \le C_{N,k}\le2Nk.
\end{equation}
Moreover,
\begin{equation}\label{eq:blockrelative}
 0\le1-\frac{C_{N,k}}{2Nk}\le6\sqrt{\frac kN}.
\end{equation}
Consequently, $C_{N,k_N}\sim2Nk_N$ whenever $1\le k_N=o(N)$.
\end{proposition}
\begin{proof}
\emph{Step 1: split the energy into block energies.}
Permuting the coordinates does not change either quadratic form. We may therefore partition them into $k$ blocks, with one positive coordinate in each block, $k-r$ blocks of size $n$, and $r$ blocks of size $n+1$. For a point $x$, denote by $x^{(i)}$ its coordinates in block $i$, and put
\[
 s_i(x)=|x^{(i)}|^2\ge1,\qquad
 c_i=C_{n_i,1},\qquad n_i\in\{n,n+1\}.
\]
The sum of all $s_i(x)$ is $|x|^2$. Let $\E_i[u]$ be the energy obtained by summing only the coordinate differences in block $i$, including its Dirichlet boundary terms. Then
$\E[u]=\sum_i\E_i[u]$.

Fix the coordinates outside block $i$. The resulting section of a finitely supported $u$ is a finitely supported function on $\N\times\Z^{n_i-1}$. Its Hardy inequality, summed over all choices of the other coordinates, yields
\[
 \E_i[u]\ge \sum_x\frac{c_i}{s_i(x)}|u(x)|^2.
\]
All sums are finite, so their order can be changed directly.

\emph{Step 2: combine the block weights.}
At each $x$, Cauchy--Schwarz applied to the two vectors
$(\sqrt{c_i/s_i(x)})_i$ and $(\sqrt{s_i(x)})_i$ gives
\[
 \left(\sum_i\sqrt{c_i}\right)^2
 \le\left(\sum_i\frac{c_i}{s_i(x)}\right)|x|^2.
\]
This remains valid if some $c_i$ is zero. Summing the block inequalities therefore gives
\[
 \E[u]\ge\sum_x|u(x)|^2\sum_i\frac{c_i}{s_i(x)}
 \ge\left(\sum_i\sqrt{c_i}\right)^2\D[u].
\]
The lower bound in \eqref{eq:blocks} follows on taking the Hardy infimum. The upper bound is Proposition~\ref{prop:fold}.

\emph{Step 3: quantify the loss from the block size.}
First assume $n\ge32$, and put $d=2\sqrt2/\sqrt n\le1/2$.
For $j=n,n+1$, equation \eqref{eq:onefold} gives
\[
 C_{j,1}\ge2j(1-d),
\]
because
\[
 \frac{\max\{2\sqrt{6j},32\}}{2j}
 =\max\left\{\frac{\sqrt6}{\sqrt j},\frac{16}{j}\right\}
 \le\frac{2\sqrt2}{\sqrt n}.
\]
For $T=(k-r)\sqrt n+r\sqrt{n+1}$, one has
\[
 T\ge k\sqrt n,\qquad N<k(n+1),\qquad
 \frac{T^2}{Nk}\ge\frac n{n+1}.
\]
Thus \eqref{eq:blocks} implies
\[
 1-\frac{C_{N,k}}{2Nk}
 \le1-(1-d)\frac n{n+1}
 \le d+\frac1{n+1}
 <\frac4{\sqrt n}.
\]
Since $n=\lfloor N/k\rfloor\ge N/(2k)$, this is less than
$4\sqrt{2k/N}\le6\sqrt{k/N}$. If $n<32$, then $N/k<32$, so $6\sqrt{k/N}>1$. In that case \eqref{eq:blockrelative} follows from $0\le C_{N,k}\le2Nk$. This proves the bound for all block sizes, and its limit proves the last assertion.
\end{proof}

\subsection{A uniform relative estimate including the puncture}
The block estimate controls small $k/N$, while the annular estimate controls large $k$. Combining them gives a relative error uniform over all sectors.

\begin{proposition}\label{prop:uniformrelative}
Put
\begin{equation}\label{eq:leadingall}
 B_{N,k}=N^2\Psi(k/N)+2N\mathbf1_{\{k=0\}},
 \qquad0\le k\le N.
\end{equation}
For every integer $N\ge1$,
\begin{equation}\label{eq:uniformrate}
 0\le1-\frac{C_{N,k}}{B_{N,k}}\le20N^{-1/3}
 \qquad(0\le k\le N).
\end{equation}

\end{proposition}
\begin{proof}
\emph{Step 1: obtain a linear lower bound for the curve.}
Fix $0<\alpha\le1/32$, let $k_N=\lfloor\alpha N\rfloor$, and restrict to $N$ large enough that $k_N\ge1$. Then $n=\lfloor N/k_N\rfloor\ge32$, so \eqref{eq:onefold} gives $C_{n,1}\ge n$ and $C_{n+1,1}\ge n+1$. The block bound yields
\[
 \frac{C_{N,k_N}}{N^2}\ge\frac{k_N}{N}\frac n{n+1}
 \ge\frac{32}{33}\frac{k_N}{N}.
\]
Proposition~\ref{prop:localization}, which already proves convergence for fixed positive proportion, now gives
$\Psi(\alpha)\ge(32/33)\alpha$. For $1/32\le\alpha\le1$, monotonicity gives
$\Psi(\alpha)\ge\Psi(1/32)\ge1/33\ge\alpha/33$.
Hence the convenient global bound
\begin{equation}\label{eq:linearpsilower}
 \Psi(\alpha)\ge\frac{\alpha}{33}\qquad(0\le\alpha\le1)
\end{equation}
holds throughout the interval.

\emph{Step 2: combine the two errors for $k\ge1$.}
Let $\rho=k/N$. The upper bounds $C_{N,k}\le N^2\Psi(\rho)\le2Nk$ and \eqref{eq:blockrelative} imply
\[
 0\le1-\frac{C_{N,k}}{N^2\Psi(\rho)}
 \le1-\frac{C_{N,k}}{2Nk}\le6\sqrt{\frac kN}.
\]
On the other hand, \eqref{eq:rate} and \eqref{eq:linearpsilower} give
\[
 1-\frac{C_{N,k}}{N^2\Psi(\rho)}
 \le\frac{4\pi}{N}\sqrt{\frac{4+\rho^{-1}}{\Psi(\rho)}}
 \le\frac{4\pi\sqrt{165}}{k},
\]
where $4+\rho^{-1}\le5/\rho$ was used. For positive $a,b$, $\min\{a,b\}\le a^{2/3}b^{1/3}$. Applying this to the last two bounds gives
\[
 1-\frac{C_{N,k}}{N^2\Psi(\rho)}
 \le(144\pi\sqrt{165})^{1/3}N^{-1/3}<20N^{-1/3}.
\]
For the last elementary comparison, $\pi<22/7$ and $\sqrt{165}<13$ give $144\pi\sqrt{165}<8000$.

\emph{Step 3: include $k=0$.}
Proposition~\ref{prop:punctured} gives
\[
 0\le1-\frac{C_{N,0}}{2N}
 \le\max\{\sqrt6\,N^{-1/2},16N^{-1}\}
 \le16N^{-1/3}.
\]
Together with Step~2 this proves \eqref{eq:uniformrate}. Since $B_{N,k}>0$ in every sector, the asserted uniform relative convergence follows.
\end{proof}

Since $20N^{-1/3}\to0$, Proposition~\ref{prop:uniformrelative} proves Theorem~\ref{thm:main}. The coefficient $\Psi(0)=0$ describes the quadratic scale; the separate term $2N\mathbf1_{\{k=0\}}$ accounts for the first-order cost of the puncture. The block estimate and the fixed-proportion limit give the regimes in Table~\ref{tab:scales}. The next section determines the corrections to $2Nk$ when $k/N\to0$.

\section{The endpoint of vanishing boundary proportion}\label{sec:smallalpha}

As $\alpha$ tends to zero, all minimizing Mathieu parameters approach zero on the scale $q\sim2\alpha$. This localization places every global minimum in the neighborhood where the implicit function theorem gives a unique analytic branch.

Put
\begin{equation}\label{eq:FandH}
 \begin{split}
 B_\alpha(q)&=(1-\alpha)a_0(q)+\alpha b_2(q),\\
 F_\alpha(q)&=\frac{[2q+B_\alpha(q)]^2}{16q},\\
 H_\alpha(q)&=2q+2qB_\alpha'(q)-B_\alpha(q).
 \end{split}
\end{equation}
Differentiation gives
\begin{equation}\label{eq:Fprime}
 F_\alpha'(q)=\frac{[2q+B_\alpha(q)]H_\alpha(q)}{16q^2}.
\end{equation}
The first bracket is positive by \eqref{eq:mathetas}. Therefore the stationary parameters are exactly the zeros of $H_\alpha$.

\begin{lemma}\label{lem:scaledmin}
For any choice of global minimizer $q_\alpha$ of $F_\alpha$,
\[
 q_\alpha\longrightarrow0,\qquad
 \frac{q_\alpha}{\alpha}\longrightarrow2
 \quad(\alpha\downarrow0).
\]
Both conclusions hold uniformly over all global minimizers. Moreover,
$\Psi(\alpha)/(2\alpha)\to1$.
\end{lemma}
\begin{proof}
\emph{Step 1: exclude every interval bounded away from zero.}
Since $b_2(q)>a_0(q)$ and $2q+a_0(q)>0$,
$F_\alpha(q)\ge F_0(q)>0$ for $q>0$. Formula \eqref{eq:smalleta}, with $\lambda=4/q$, gives $F_0(q)\to1/4$ as $q\to\infty$. Continuity and positivity on compact subintervals then imply
\[
 \inf_{q\ge\eps}F_0(q)>0\qquad\hbox{for each }\eps>0.
\]
On the other hand, every global minimum has value $\Psi(\alpha)<2\alpha$. Consequently no sequence of global minimizers can remain in $[\eps,\infty)$ as $\alpha\downarrow0$. This proves convergence to zero for all choices of minimizer.

\emph{Step 2: identify a compact range after rescaling.}
The small-$q$ expansions, derived in Appendix~\ref{app:coefficients}, are
\begin{equation}\label{eq:mathieusmall}
 a_0(q)=-\frac{q^2}{2}+\frac{7q^4}{128}+O(q^6),\qquad
 b_2(q)=4-\frac{q^2}{12}+\frac{5q^4}{13824}+O(q^6).
\end{equation}
In particular, for some $q_0>0$,
\[
 2q+a_0(q)\ge q,\qquad b_2(q)-a_0(q)\ge2
 \quad(0<q\le q_0).
\]
By Step 1, every minimizing parameter lies in this interval for all sufficiently small $\alpha$. At such a parameter,
\[
 \frac{(q_\alpha+2\alpha)^2}{16q_\alpha}
 \le F_\alpha(q_\alpha)<2\alpha.
\]
The separate nonnegative terms $q_\alpha/16$ and $\alpha^2/(4q_\alpha)$ in the left side give the convenient bounds
\[
 \frac{\alpha}{8}\le q_\alpha\le32\alpha.
\]
These bounds give a fixed compact interval for the rescaled parameter.

\emph{Step 3: pass to the rescaled variational problem.}
For $z$ in the fixed compact interval $[1/8,32]$, equation \eqref{eq:mathieusmall} gives, uniformly in $z$,
\begin{equation}\label{eq:scaledlimit}
 \frac{F_\alpha(\alpha z)}{\alpha}
 \longrightarrow \frac{(2z+4)^2}{16z}
 =:\Phi(z)=\frac z4+1+\frac1z.
\end{equation}
The identity $\Phi(z)-2=(z-2)^2/(4z)$ proves that its minimum is uniquely attained at $z=2$.

To check that this identifies all minima, take any sequence $\alpha_j\downarrow0$ and corresponding global minima. A subsequence of $z_j=q_{\alpha_j}/\alpha_j$ converges to $z_0\in[1/8,32]$. The minimizing inequality
$F_{\alpha_j}(\alpha_jz_j)\le F_{\alpha_j}(2\alpha_j)$ and uniform convergence in \eqref{eq:scaledlimit} imply $\Phi(z_0)\le\Phi(2)$. Hence $z_0=2$. Every convergent subsequence has the same limit, proving $q_\alpha/\alpha\to2$ uniformly over the minima. Substitution in \eqref{eq:scaledlimit} also gives $\Psi(\alpha)/\alpha\to\Phi(2)=2$.
\end{proof}

\begin{proposition}\label{prop:smallalpha}
For all sufficiently small $\alpha>0$, $F_\alpha$ has a unique global minimizing parameter $q_\alpha$. Both $q_\alpha$ and $\Psi(\alpha)$ extend real analytically across $\alpha=0$, with
\begin{align}
 q_\alpha&=2\alpha+3\alpha^2+\frac{13}{2}\alpha^3
                 +\frac{251}{16}\alpha^4+O(\alpha^5),\label{eq:qexp}\\
 \Psi(\alpha)&=2\alpha-\alpha^2-\frac{\alpha^3}{6}
                 -\frac{7\alpha^4}{48}+O(\alpha^5).\label{eq:psiexp}
\end{align}
In addition,
\begin{equation}\label{eq:lambdaexp}
 \lambda_\alpha=4/q_\alpha
 =\frac2\alpha-3-2\alpha-\frac{47}{16}\alpha^2+O(\alpha^3).
\end{equation}
\end{proposition}
\begin{proof}
\emph{Step 1: construct the unique branch containing all minima.}
The characteristic values in \eqref{eq:mathieusmall} are real analytic near zero. After the substitution $q=\alpha z$, every term of $H_\alpha(\alpha z)$ contains a factor $\alpha$. Thus
\[
 K(\alpha,z):=\frac{H_\alpha(\alpha z)}{\alpha}
\]
extends real analytically to a neighborhood of $(0,2)$, with
\[
 K(0,z)=2z-4,\qquad \partial_zK(0,2)=2.
\]
The implicit function theorem gives an analytic function $z(\alpha)$, with $z(0)=2$, and a neighborhood in which $K(\alpha,z)=0$ if and only if $z=z(\alpha)$. Lemma~\ref{lem:scaledmin} places every global minimizer in this neighborhood for small positive $\alpha$. Each minimizer is interior and hence stationary by \eqref{eq:Fprime}; therefore every one equals $q_\alpha=\alpha z(\alpha)$. Existence was already proved in Lemma~\ref{lem:attainment}. This establishes global uniqueness and the analytic extension of $q_\alpha$.

\emph{Step 2: compute its coefficients.}
Equation \eqref{eq:mathieusmall} gives
\begin{equation}\label{eq:Hexp}
 H_\alpha(q)=2q-4\alpha+
 \left(-\frac32+\frac54\alpha\right)q^2
 +\frac{49}{128}q^4+O(\alpha q^4+q^6).
\end{equation}
Here the remainder is uniform for $(\alpha,q)$ near $(0,0)$. Along $q=O(\alpha)$ it is $O(\alpha^5)$, which suffices to determine the four displayed coefficients of $q_\alpha$.

Write
$q_\alpha=2\alpha+A\alpha^2+B\alpha^3+C\alpha^4+O(\alpha^5)$.
The contributions through order four are
\[
 q_\alpha^2=4\alpha^2+4A\alpha^3+(A^2+4B)\alpha^4+O(\alpha^5),
 \qquad q_\alpha^4=16\alpha^4+O(\alpha^5).
\]
The coefficients of $\alpha^2,\alpha^3,\alpha^4$ in $H_\alpha(q_\alpha)=0$ are respectively
\[
 2A-6,\qquad 2B-6A+5,\qquad
 2C-\frac32(A^2+4B)+5A+\frac{49}{8}.
\]
Solving in that order gives $A=3$, $B=13/2$, and $C=251/16$, proving \eqref{eq:qexp}.

\emph{Step 3: evaluate the minimizing value.}
Substitution into $2q+B_\alpha(q)$ gives
\[
 2q_\alpha+B_\alpha(q_\alpha)
 =8\alpha+4\alpha^2+\frac{26}{3}\alpha^3
   +\frac{79}{4}\alpha^4+O(\alpha^5).
\]
Its square divided by $16q_\alpha$ yields \eqref{eq:psiexp}; expanding $4/q_\alpha$ yields \eqref{eq:lambdaexp}. Analyticity follows from the factorization
$q_\alpha=\alpha z(\alpha)$ with $z(0)=2$ and
$2q_\alpha+B_\alpha(q_\alpha)=\alpha w(\alpha)$ with analytic $w$ and $w(0)=8$. Hence
\[
 \Psi(\alpha)=\alpha\,\frac{w(\alpha)^2}{16z(\alpha)}
\]
has an analytic extension across zero. This also justifies all divisions used in the coefficient calculation.
\end{proof}

The derivative $\Psi'(0)=2$ explains the agreement between the curve and the block law $2Nk$. The uniform relative estimate in Section~\ref{sec:blocks} has already been proved independently. The remaining coefficients describe the lattice constants when the correction terms exceed the order-$N$ error from radial localization.

\begin{corollary}\label{cor:sublinearcorrections}
For every sequence $1\le k=k_N=o(N)$,
\begin{equation}\label{eq:sublinearexpansion}
 C_{N,k}=2Nk-k^2-\frac{k^3}{6N}-\frac{7k^4}{48N^2}
          +O\!\left(N+\frac{k^5}{N^3}\right).
\end{equation}
More precisely, each of the following assertions holds under its stated hypotheses:
\begin{align}
 C_{N,k}&=2Nk-k^2+o(k^2),
 &&\sqrt N=o(k),\quad k=o(N),\label{eq:twotermsublinear}\\
 C_{N,k}&=2Nk-k^2-\frac{k^3}{6N}+o(k^3/N),
 &&N^{2/3}=o(k),\quad k=o(N),\label{eq:threetermsublinear}\\
 C_{N,k}&=2Nk-k^2-\frac{k^3}{6N}-\frac{7k^4}{48N^2}+o(k^4/N^2),
 &&N^{3/4}=o(k),\quad k=o(N).\label{eq:fourtermsublinear}
\end{align}
\end{corollary}
\begin{proof}
Analyticity of $\Psi$ at zero gives constants $\alpha_0>0$ and $A<\infty$ such that the remainder after the fourth-degree polynomial in \eqref{eq:psiexp} has absolute value at most $A\alpha^5$ for $0\le\alpha\le\alpha_0$. Write
\[
 C_{N,k}=N^2\Psi(k/N)-R_{N,k},\qquad
 0\le R_{N,k}\le4\pi\sqrt{10}\,N
\]
by \eqref{eq:rate}. Multiplication of the Taylor expansion by $N^2$ now proves \eqref{eq:sublinearexpansion}, with a uniform implied constant whenever $k/N\le\alpha_0$.

For \eqref{eq:twotermsublinear}, truncate the same expansion after $\alpha^2$. The remainder is $O(N+k^3/N)$; divided by $k^2$, it is bounded by a constant times $N/k^2+k/N$, which tends to zero under the stated hypotheses. For \eqref{eq:threetermsublinear}, the remainder is $O(N+k^4/N^2)$. Dividing by $k^3/N$ gives the bound $O(N^2/k^3+k/N)$, again tending to zero. Finally, divide the remainder in \eqref{eq:sublinearexpansion} by $k^4/N^2$. The resulting bound is $O(N^3/k^4+k/N)$, which proves \eqref{eq:fourtermsublinear}.
\end{proof}

The growth conditions place each displayed correction above the localization error. For fixed $k$, the block argument gives the leading law; a constant-order correction remains open.

The oscillator--Mathieu identities, Lemma~\ref{lem:attainment}, Proposition~\ref{prop:shape}, and Proposition~\ref{prop:smallalpha} together prove Theorem~\ref{thm:curve}.

\section{The endpoint energy--moment profile}\label{sec:kappa}

At $\alpha=1$, uniqueness selects one normalized positive minimizer of the one-dimensional energy--moment product. Nondegeneracy then yields the analytic branch near one in Corollary~\ref{cor:nearone}. Tensor powers of this profile give the asymptotically optimal Hardy trials constructed in Section~\ref{sec:localization}. The proof has three steps. Interpolation and an oscillator comparison place every minimizing parameter in $20/9<q<4$. On $0<q\le4$, a spectral gap greater than $12$ and an adjacency variance less than $1$ give $b_2''(q)>-1/6$. This curvature bound makes the stationary function strictly increasing on the interval containing all minima and also proves nondegeneracy.

For $f\in D(X)\subset\ell^2(\N)$, write $e(f)=\langle f,L_+f\rangle$ and $m(f)=\|Xf\|_2^2$.

\Needspace{17\baselineskip}
\begin{theorem}[Exact characterization of the endpoint]\label{thm:kappacharacterization}
One has
\[
 \kappa=\min_{\substack{f\in D(X)\\\|f\|_2=1}}e(f)m(f)
       =\min_{q>0}\frac{[2q+b_2(q)]^2}{16q}.
\]
The positive normalized minimizing profile $f_*$ is unique, and every complex minimizer is a constant phase multiple of $f_*$. The minimizing parameter $q_*$ is nondegenerate and is the unique zero in $(0,4]$ of
\begin{equation}\label{eq:root}
 2q+2qb_2'(q)-b_2(q)=0.
\end{equation}
It satisfies
\begin{equation}\label{eq:qinterval}
 \frac{20}{9}<q_*<4,\qquad
 \kappa=\frac{[2q_*+b_2(q_*)]^2}{16q_*}
 \le\frac{105862}{59049}<\frac{4489}{2500}<\frac95.
\end{equation}
Moreover, $f_*$ is the ground state of $L_++\lambda_*X^2$ with $\lambda_*=4/q_*$, and $e(f_*)=\lambda_*m(f_*)$.
\end{theorem}

\subsection{Interpolation and bounds for minimizing profiles}
Proposition~\ref{prop:uncertainty} gives
\[
 \kappa=\inf_{\substack{f\in D(X)\\\|f\|_2=1}}e(f)m(f)
       =\min_{\lambda>0}\frac{\eta_+(\lambda)^2}{4\lambda}.
\]
A useful strict upper bound comes from the finitely supported vector
$v=(1,1/4,1/40,0,\ldots)$. Direct summation gives
\[
 \|v\|_2^2=\frac{1701}{1600},\qquad
 \ip v{L_+v}=\frac{1291}{800},\qquad
 \|Xv\|_2^2=\frac{2009}{1600}.
\]
Consequently,
\begin{equation}\label{eq:U}
 \kappa\le\frac{\ip v{L_+v}\|Xv\|_2^2}{\|v\|_2^4}
 =U:=\frac{105862}{59049}<\frac{4489}{2500}<\frac95.
\end{equation}

\begin{lemma}\label{lem:interpolation}
Every normalized $f\in D(X)\subset\ell^2(\N)$ satisfies
\begin{equation}\label{eq:eminterpolation}
 e(f)m(f)\ge\frac94-\frac{11}{12}e(f)+\frac{e(f)^2}{16}.
\end{equation}
\end{lemma}
\begin{proof}
First suppose that $f$ is finitely supported. Put $f_0=0$ and let $F$ be its piecewise linear interpolation on $[0,\infty)$. On $[j,j+1]$,
$F(j+t)=(1-t)f_j+tf_{j+1}$ for $0\le t\le1$. It follows that
\[
 \int_j^{j+1}|F'|^2=|f_{j+1}-f_j|^2,\qquad
 \int_j^{j+1}|F|^2=
 \frac{|f_j|^2+|f_{j+1}|^2+\Rea(f_j\overline{f_{j+1}})}3.
\]
Summing, and using $e=2-2\Rea\sum_{j\ge1}f_j\overline{f_{j+1}}$, gives
\begin{equation}\label{eq:interp}
 \int_0^\infty|F'|^2=e,\qquad
 \int_0^\infty|F|^2=1-\frac e6,\qquad
 \int_0^\infty x^2|F|^2\le m+\frac16.
\end{equation}
To justify the last estimate, convexity gives
$|F(j+t)|^2\le(1-t)|f_j|^2+t|f_{j+1}|^2$.
The coefficient of $|f_j|^2$ after summing the two adjacent intervals is
\[
 \int_{-1}^1(j+t)^2(1-|t|)\,dt=j^2+\frac16,
\]
also for $j=1$. The coefficient at $j=0$ contributes nothing because $f_0=0$.

We next use a continuous half-line uncertainty inequality, keeping its constant explicit. For $b>0$, expansion and integration by parts give
\[
 \int_0^\infty\left|F'-\left(\frac1x-bx\right)F\right|^2
 =\int_0^\infty|F'|^2+b^2\int_0^\infty x^2|F|^2
       -3b\int_0^\infty|F|^2.
\]
There is no boundary term at infinity because $F$ has compact support. Near zero $F(x)=xf_1$, so $(1/x-bx)|F(x)|^2\to0$; thus the integration by parts at zero is justified as well. If $E_c,M_c,A_c$ denote the three integrals on the right, nonnegativity of the left side for all $b>0$, evaluated at $b=3A_c/(2M_c)$, gives
\[
 E_cM_c\ge\frac94 A_c^2.
\]
Combining this inequality with \eqref{eq:interp} yields
$e(m+1/6)\ge(9/4)(1-e/6)^2$. Expanding and subtracting $e/6$ gives \eqref{eq:eminterpolation}.

For a general $f\in D(X)$, its finite truncations converge in the norm
$\|f\|_2+\|Xf\|_2$. Normalize these truncations. Their energies converge because $L_+$ is bounded, and their second moments converge by construction. Passing to the limit proves the claimed inequality for $f$.
\end{proof}

\begin{lemma}\label{lem:roughlocalization}
The product $e(f)m(f)$ attains its infimum on the unit sphere of $D(X)$. Every minimizing profile is a ground state of $L_++\lambda_*X^2$ with $\lambda_*=e(f)/m(f)$. Every minimizing oscillator parameter satisfies
\begin{equation}\label{eq:roughlambda}
 \frac{81}{605}<\lambda_*<\frac95.
\end{equation}
\end{lemma}
\begin{proof}
Lemma~\ref{lem:attainment} gives a minimizing oscillator parameter, and Lemma~\ref{lem:balance} gives a normalized ground profile with product $em=\kappa$. Thus the profile infimum is attained.

Conversely, let $f$ be any normalized minimizing profile and set $\lambda_*=e/m$. Both quantities are positive. The variational principle and the definition of $\kappa$ give
\[
 \kappa\le\frac{\eta_+(\lambda_*)^2}{4\lambda_*}
 \le\frac{(e+\lambda_*m)^2}{4\lambda_*}
 =em=\kappa.
\]
Equality holds throughout. In particular, $f$ attains the ground-state Rayleigh quotient. It belongs to the lowest eigenspace, hence to $D(X^2)$, and satisfies
\begin{equation}\label{eq:euler}
 (L_++\lambda_*X^2)f=2ef,\qquad \lambda_*=\frac em.
\end{equation}
This also proves that its associated parameter minimizes the scalar problem. Conversely, the ground profile at any minimizing scalar parameter satisfies $em=\kappa$ and $\lambda_*=e/m$ by Lemma~\ref{lem:balance}.

It remains to bound these parameters. Lemma~\ref{lem:interpolation} and $em=\kappa<9/5$ imply
\[
 \frac94-\frac{11}{12}e<\frac95,\qquad
 e>\frac{27}{55},\qquad
 m<\frac{9/5}{27/55}=\frac{11}{3}.
\]
Therefore $\lambda_*>(27/55)/(11/3)=81/605$. Since $m\ge1$, one also has $\lambda_*=em/m^2\le\kappa<9/5$. These bounds apply to every minimizing profile and every minimizing scalar parameter.
\end{proof}

\subsection{An analytic interval containing all minimizing parameters}
The next comparison raises the lower bound on $\lambda_*$ to $1$, and hence places $q_*=4/\lambda_*$ below $4$. The coefficient $67/25$ is chosen so that $(67/25)^2/4>U$.

\begin{lemma}\label{lem:comparison}
If $c=67/25$ and $1/3\le s\le1$, then
\begin{equation}\label{eq:comparison}
 L_++s^2X^2-csI>0.
\end{equation}
Here strict positivity means that the operator is bounded below by a positive multiple of the identity for each such $s$.
\end{lemma}
The proof in Appendix~\ref{app:comparison} separates the first three coordinates from the tail and reduces the remaining sign to the concavity of an explicit polynomial.

\begin{corollary}\label{cor:localization}
Every minimizing parameter satisfies
\begin{equation}\label{eq:loclambda}
 1<\lambda_*<\frac95,\qquad
 \frac{20}{9}<q_*:=4/\lambda_*<4.
\end{equation}
\end{corollary}
\begin{proof}
Suppose $\lambda_*\le1$. Since $81/605>1/9$, Lemma~\ref{lem:roughlocalization} gives $1/3<\sqrt{\lambda_*}\le1$. Apply Lemma~\ref{lem:comparison} with $s=\sqrt{\lambda_*}$. The variational principle yields $\eta_+(\lambda_*)>c\sqrt{\lambda_*}$, and therefore
\[
 \kappa=\frac{\eta_+(\lambda_*)^2}{4\lambda_*}
 >\frac{c^2}{4}=\frac{4489}{2500}>U.
\]
This contradicts \eqref{eq:U}. The remaining bounds follow from \eqref{eq:roughlambda} and the decreasing change of variables $q=4/\lambda$.
\end{proof}

\subsection{The gap and the ground-state tail}
The sine expansion of \eqref{eq:mathieuop} on $(0,\pi/2)$, followed by alternating the signs of the coefficients, identifies $b_2(q)$ with the ground eigenvalue of
\begin{equation}\label{eq:Tq}
 T_q=D-qA,\qquad
 D=\operatorname{diag}(4,16,36,\ldots),\qquad A=S+S^*.
\end{equation}
Here $S$ is the unilateral shift on $\ell^2(\N)$ and $\|A\|\le2$. The domain of $T_q$ is $D(D)$, independent of $q$. Its lowest eigenvalue is simple for $q>0$, with a positive normalized eigenvector $f$. In this subsection write $b=b_2(q)$.

For $q>0$, the trial vector $\delta_1+t\delta_2$ has quotient below $4$ when $t>0$ is sufficiently small: after subtracting $4(1+t^2)$, its numerator is $12t^2-2qt$. Thus $b<4$. At $q=0$, direct perturbation of the simple eigenvalue $4$ gives $b_2(0)=4$ and $b_2'(0)=-\ip{\delta_1}{A\delta_1}=0$.

\begin{lemma}\label{lem:gap}
For $0<q\le4$, the second eigenvalue of $T_q$ is at least $16$. Hence its ground-state gap $g$ satisfies
\begin{equation}\label{eq:gap}
 g\ge16-b>12.
\end{equation}
\end{lemma}
\begin{proof}
Split a vector into its first coordinate $x$, its second coordinate $y$, and its tail $z$ on $\{3,4,\ldots\}$. Let $C$ be the restriction of $T_q-16I$ to that tail with Dirichlet boundary at coordinate two. Since the diagonal starts at $20$,
\[
 C\ge(20-2q)I\ge12I.
\]
Thus $C^{-1}$ exists and $\ip{\delta_3}{C^{-1}\delta_3}\le1/12$.

For vectors in the form domain, completion of squares gives the exact identity
\[
 \begin{split}
 \ip{(x,y,z)}{(T_q-16I)(x,y,z)}
 ={}&-12\left|x+\frac q{12}y\right|^2\\
 &+\ip{z-qyC^{-1}\delta_3}{C(z-qyC^{-1}\delta_3)}\\
 &+q^2\left(\frac1{12}-\ip{\delta_3}{C^{-1}\delta_3}\right)|y|^2.
 \end{split}
\]
The identity is a form identity; its middle term denotes the square of the $C^{1/2}$ norm. The shift by $C^{-1}\delta_3$ is legitimate because this vector lies in $D(C)$. Both terms after the first are nonnegative. On the codimension-one subspace $x+qy/12=0$, the full form is therefore nonnegative. By the min--max principle, $T_q-16I$ has at most one negative eigenvalue, counting multiplicity. Compact resolvent then gives the asserted bound on the second eigenvalue. Subtracting $b<4$ proves \eqref{eq:gap}.
\end{proof}

\begin{lemma}\label{lem:tail}
Let $0<q\le4$ and set $r_j=f_{j+1}/f_j$. Then
\begin{equation}\label{eq:ratios}
 0<r_1<\frac12,\qquad
 0<r_j\le\frac{r_1}{2}<\frac14\quad(j\ge2).
\end{equation}
If $p=f_1^2$, then
\begin{equation}\label{eq:pbound}
 1-p\le\frac{16}{15}r_1^2p,\qquad p\ge\frac{15}{19}.
\end{equation}
Moreover,
\begin{equation}\label{eq:variance}
 \Var_f(A):=\ip f{A^2f}-\ip f{Af}^2
 \le1-\frac{2309}{5415}r_1^2<1.
\end{equation}
\end{lemma}
\begin{proof}
\emph{Step 1: bound ratios using a tail resolvent.}
For $j\ge1$, restrict the eigenvalue equation to coordinates $j+1,j+2,\ldots$, and denote that restriction of $T_q$ by $T^{(j+1)}_q$. Then
\[
 (T^{(j+1)}_q-b)f_{\ge j+1}=qf_j\delta_{j+1}.
\]
The operator on the left is bounded below by
$4(j+1)^2-2q-b>0$, because $j\ge1$, $q\le4$, and $b<4$. Inverting and taking its first coordinate gives
\begin{equation}\label{eq:tailratio}
 r_j=q\ip{\delta_{j+1}}{(T^{(j+1)}_q-b)^{-1}\delta_{j+1}}
 \le\frac{q}{4(j+1)^2-2q-b}.
\end{equation}
All ratios are positive since $f>0$. In particular $r_j\le1/6$ for $j\ge2$.

\emph{Step 2: compare the whole tail with the first ratio.}
The first recurrence gives $d:=4-b=qr_1>0$. The second recurrence gives
\[
 r_1=\frac{q}{12+d-qr_2}
 \le\frac{q}{12-4/6}=\frac{3q}{34}
 \le\frac6{17}<\frac12.
\]
It follows that $d=qr_1<2$. The same recurrence also gives
$r_1\ge q/(12+d)$. From the third recurrence and the preliminary estimate $r_3\le1/6$,
\[
 r_2=\frac{q}{32+d-qr_3}\le\frac{3q}{94},\qquad
 \frac{r_2}{r_1}\le\frac{3(12+d)}{94}<\frac{21}{47}<\frac12.
\]
For $j\ge3$, \eqref{eq:tailratio} gives $r_j\le q/52$, so
\[
 \frac{r_j}{r_1}\le\frac{12+d}{52}<\frac{14}{52}<\frac12.
\]
These estimates prove \eqref{eq:ratios}.

\emph{Step 3: sum the tail and bound the variance.}
Since $f_2=r_1f_1$ and every subsequent ratio is at most $1/4$,
\[
 1-p=\sum_{j\ge2}f_j^2
 \le p r_1^2\sum_{\ell\ge0}16^{-\ell}
 =\frac{16}{15}p r_1^2.
\]
Thus $p\ge[1+(16/15)r_1^2]^{-1}\ge15/19$, proving \eqref{eq:pbound}. Also $r_jr_{j+1}\le r_1^2/2$ for every $j\ge1$, and hence
\[
 \sum_{j\ge1}f_jf_{j+2}
 \le\frac{r_1^2}{2}\sum_{j\ge1}f_j^2=\frac{r_1^2}{2}.
\]
The diagonal of $A^2$ is $1$ at the first coordinate and $2$ thereafter. Its remaining entries connect coordinates two steps apart. Therefore
\[
 \ip f{A^2f}=2-p+2\sum_{j\ge1}f_jf_{j+2},
 \qquad \ip f{Af}=2\sum_{j\ge1}f_jf_{j+1}\ge2r_1p.
\]
Using $1-p\le(16/15)r_1^2$ and $p\ge15/19$, we obtain
\[
 \Var_f(A)\le1+\frac{31}{15}r_1^2
             -4r_1^2\left(\frac{15}{19}\right)^2
 =1-\frac{2309}{5415}r_1^2<1.
\]
This proves \eqref{eq:variance}.
\end{proof}

\subsection{Curvature, uniqueness, and nondegeneracy}
The gap controls the reduced resolvent in the second derivative of $b_2$. Combining it with the variance bound gives the curvature estimate needed for the stationary equation.
\begin{lemma}\label{lem:curvature}
For $0<q\le4$,
\begin{equation}\label{eq:curvature}
 b_2''(q)>-\frac16,\qquad b_2'(q)>-\frac q6.
\end{equation}
\end{lemma}
\begin{proof}
The bounded perturbation $-qA$ of $D$ is an analytic operator family with fixed domain. Simplicity allows a locally analytic, real normalized eigenvector $f(q)$; see \cite[Chapter~VII]{Kato}. Differentiating its normalization gives $\ip f{f'}=0$. Differentiation of $T_qf=bf$ and then scalar multiplication by $f$ give
\[
 b'=-\ip f{Af},\qquad
 (T_q-b)f'=(A-\ip f{Af})f.
\]
Let $v=(A-\ip f{Af})f$. This vector is orthogonal to $f$. On $f^\perp$, the reduced inverse
$R=(T_q-b)^{-1}_{f^\perp}$ is positive and has norm $1/g$. Thus $f'=Rv$, and a second differentiation yields
\[
 b''=-2\ip{f'}{Af}
     =-2\ip v{Rv}
     \ge-\frac{2\|v\|_2^2}{g}
     =-\frac{2\Var_f(A)}g.
\]
Lemmas~\ref{lem:gap} and \ref{lem:tail} give $g>12$ and $\Var_f(A)<1$, so $b''>-1/6$. The eigenvalue is analytic also at $q=0$, where $b_2'(0)=0$. Integration of the strict curvature bound on $(0,q)$ gives $b_2'(q)>-q/6$.
\end{proof}

\begin{proof}[Proof of Theorem~\ref{thm:kappacharacterization}]
Set
\[
 F(q)=\frac{[2q+b_2(q)]^2}{16q},\qquad
 H(q)=2q+2qb_2'(q)-b_2(q).
\]
Lemma~\ref{lem:curvature} implies
\begin{equation}\label{eq:Hmonotone}
 H'(q)=2+b_2'(q)+2qb_2''(q)>2-\frac q2\ge0
 \qquad(0<q\le4).
\end{equation}
In particular $H$ is strictly increasing on this whole interval, including at its right endpoint.

Existence of at least one global minimum follows from Lemma~\ref{lem:attainment}. Corollary~\ref{cor:localization} places every global minimum in $(20/9,4)$, and \eqref{eq:Fprime} shows that each is a zero of $H$. A strictly increasing function has at most one zero. Hence there is exactly one global minimizing parameter $q_*$, and it is also the unique zero of $H$ in $(0,4]$. The localization supplies the global conclusion.

Differentiating \eqref{eq:Fprime} at $q_*$ eliminates all terms containing $H(q_*)=0$, leaving
\[
 F''(q_*)=\frac{[2q_*+b_2(q_*)]H'(q_*)}{16q_*^2}>0.
\]
The first factor is positive by \eqref{eq:mathetas}, and the derivative is positive by \eqref{eq:Hmonotone}. This proves nondegeneracy.

Every nonnegative minimizing profile solves \eqref{eq:euler} at the parameter $\lambda_*=4/q_*$. Simplicity of the oscillator ground eigenvalue gives a unique positive normalized profile $f_*$. For a complex minimizing profile $f$, replacing it by $|f|$ preserves its norm and moment and cannot increase its energy. Since its moment is positive and the original product is minimal, equality in energy must hold. The energy difference is the sum of the nonnegative edge terms
\[
 2\sum_{j\ge1}\bigl(|f_j||f_{j+1}|-\Rea(\overline f_jf_{j+1})\bigr).
\]
Its vanishing forces each term to vanish. Since $|f|=f_*>0$ at every coordinate, adjacent values have the same phase. Connectedness of the half-line gives one constant phase for the entire sequence. Finally, \eqref{eq:U} and Corollary~\ref{cor:localization} supply the exact bounds in \eqref{eq:qinterval}.
\end{proof}

\subsection{Continuation near the fully folded endpoint}
\begin{corollary}\label{cor:nearone}
There is $\eps>0$ such that, for $1-\eps<\alpha\le1$, the minimum in \eqref{eq:mathieupsi} is unique. Its parameter and value extend real analytically across $\alpha=1$. In particular,
\begin{equation}\label{eq:derivativeone}
 \Psi'(1^-)=
 \frac{[2q_*+b_2(q_*)][b_2(q_*)-a_0(q_*)]}{8q_*}>0.
\end{equation}
\end{corollary}
\begin{proof}
Lemma~\ref{lem:attainment} confines all minimizing $\lambda$ parameters to a fixed compact subinterval of $(0,\infty)$ when $\alpha$ is close to one. The substitution $q=4/\lambda$ gives the same assertion for the $q$ parameters.

Take any sequence $\alpha_j\to1$ with $\alpha_j\le1$ and any choice of global minimizers $q_j$. Every convergent subsequence has a limit $\bar q>0$. The minimizing inequalities
$F_{\alpha_j}(q_j)\le F_{\alpha_j}(q)$, valid for each fixed $q>0$, pass to the limit by continuity and give $F_1(\bar q)\le F_1(q)$. Thus $\bar q=q_*$ by Theorem~\ref{thm:kappacharacterization}. It follows that all global minima approach $q_*$.

The functions $a_0,b_2$ are real analytic near the positive point $q_*$. Equation \eqref{eq:Hmonotone} gives $\partial_qH_1(q_*)>0$. The implicit function theorem therefore produces a unique analytic solution of $H_\alpha(q)=0$ near $(1,q_*)$. All global minima are eventually in that neighborhood and are stationary, so they belong to this branch. This proves uniqueness and the analytic extension of the parameter and of its value. The variational problem selects this analytic branch on $\alpha\le1$.

Lastly, differentiate $\Psi(\alpha)=F_\alpha(q_\alpha)$. The term involving $q_\alpha'$ vanishes because $\partial_qF_\alpha(q_\alpha)=0$. The remaining partial derivative at $\alpha=1$ is \eqref{eq:derivativeone}, and its sign follows from $b_2(q_*)>a_0(q_*)$ and $2q_*+b_2(q_*)>0$.
\end{proof}

\section{Continuum comparison and spectral thresholds}\label{subsec:strictcomparison}

Table~\ref{tab:endpoints} places the uniform comparison beside the known endpoint classifications.

\begin{table}[htbp]
\caption{Comparison of $C_{N,k}$ with $A_{N,k}$. The symbols $=$, $<$, and $\le$ compare these two constants. An entry $\le$ marks cases for which the table leaves equality versus strictness open. The endpoint columns are known \cite{AlpayLow,LizamaLow,LizamaPositive,LizamaTransition}; the uniform mixed-sector conclusion for $N\ge9$ is Theorem~\ref{thm:comparison}. For $(N,k)=(8,1)$, Proposition~\ref{prop:strictcriterion} gives $C_{8,1}<A_{8,1}=16$.}\label{tab:endpoints}
\centering
\setlength{\tabcolsep}{12pt}
\begin{tabular}{@{}lccc@{}}
\toprule
Dimension & $k=0$ & $1\le k<N$ & $k=N$\\
\midrule
$N=1$ & $=$ & --- & $=$\\
$N=2$ & $=$ & $\le$ & $=$\\
$3\le N\le8$ & $=$ & $\le$ & $<$\\
$N\ge9$ & $<$ & $<$ & $<$\\
\bottomrule
\end{tabular}
\end{table}

\subsection{Strict comparison}
We now prove Theorem~\ref{thm:comparison} by finitely supported trials. A two-point argument first gives a criterion that is useful in individual sectors. A product trial with a variance correction then covers every mixed sector from dimension nine onward.

\begin{proposition}[A criterion for strict comparison]\label{prop:strictcriterion}
For $1\le k\le N$, one has $C_{N,k}<2Nk$. Consequently,
\begin{equation}\label{eq:strictcriterion}
 (N-2k)^2\ge4(N+2k-1)
 \quad\Longrightarrow\quad C_{N,k}<A_{N,k}.
\end{equation}
\end{proposition}
\begin{proof}
Let $x=(1,\ldots,1,0,\ldots,0)$, with $k$ positive coordinates, and let $y=x+e_1$. Both points belong to the domain, are adjacent, and satisfy $|x|^2=k$ and $|y|^2=k+3$. For $u=\delta_x+\eps\delta_y$, direct evaluation of the two forms gives
\[
 \E[u]=2N(1+\eps^2)-2\eps,\qquad
 \D[u]=\frac1k+\frac{\eps^2}{k+3}.
\]
Hence
\[
 \E[u]-2Nk\D[u]=-2\eps+\frac{6N}{k+3}\eps^2<0
 \qquad\left(0<\eps<\frac{k+3}{3N}\right).
\]
This proves the strict bound. Finally,
\[
 4(A_{N,k}-2Nk)=(N-2k)^2-4(N+2k-1).
\]
The hypothesis in \eqref{eq:strictcriterion} makes $A_{N,k}\ge2Nk$, so the conclusion follows even when equality holds in that hypothesis.
\end{proof}

The criterion is effective away from $k\approx N/2$. To cover all sectors in dimension nine and above, we retain more information about the inverse-square moment of a product trial. The following elementary estimate measures the gain over the usual arithmetic-mean bound.

\begin{lemma}[A variance correction for the harmonic mean]\label{lem:harmonicmean}
Let $S$ be a positive random variable with finite third moment and finite $\mathbb E(S^{-1})$. Put
\[
 m=\mathbb E S,\qquad v=\mathbb E(S-m)^2,\qquad
 w=\mathbb E(S-m)^3.
\]
If $v>0$, then $mv+w>0$ and
\begin{equation}\label{eq:harmonicmean}
 \frac1{\mathbb E(S^{-1})}
 \le m-\frac{v^2}{mv+w}.
\end{equation}
\end{lemma}
\begin{proof}
The denominator is $mv+w=\mathbb E[S(S-m)^2]>0$. For any real $t$, the random variable $Q_t=1-t(S-m)$ has expectation one. Cauchy--Schwarz gives
\[
 1=(\mathbb E Q_t)^2
 \le\mathbb E(SQ_t^2)\,\mathbb E(S^{-1}).
\]
Expanding the first factor gives
\[
 \mathbb E(SQ_t^2)=m-2tv+t^2(mv+w).
\]
Its minimum over $t\in\R$ is attained at $t=v/(mv+w)$ and equals the right-hand side of \eqref{eq:harmonicmean}.
\end{proof}

\begin{proof}[Proof of Theorem~\ref{thm:comparison}]
One product family covers $N\ge10$ and, in dimension nine, $k\ge4$. Proposition~\ref{prop:strictcriterion} handles $k=1$ there; two adjusted estimates handle $k=2,3$. The known punctured-lattice comparison completes the proof for $k=0$ and establishes optimality of the uniform dimension.

\emph{Step 1: one product family and its moments.}
First suppose $k\ge1$ and put $\alpha=k/N$. On the half-line and on the line, respectively, take
\[
 f=\frac{5\delta_1+\delta_2}{\sqrt{26}},\qquad
 g=\frac{\delta_{-1}+3\delta_0+\delta_1}{\sqrt{11}},\qquad
 u=f^{\otimes k}\otimes g^{\otimes(N-k)}.
\]
The one-dimensional energies are $21/13$ and $10/11$. Under the probability distribution $|u(x)|^2$, the squared radius $S=|x|^2$ has the law
\[
 S=k+3I+J,\qquad
 I\sim\operatorname{Bin}(k,1/26),\quad
 J\sim\operatorname{Bin}(N-k,2/11),
\]
where $I,J$ are independent. In particular, $S\ge k>0$. Define the following four affine functions of $\alpha$:
\begin{align}
 e(\alpha)&=\frac{10}{11}+\frac{101}{143}\alpha,
 &m(\alpha)&=\frac2{11}+\frac{267}{286}\alpha,\notag\\
 v(\alpha)&=(1-\alpha)\frac{18}{121}
                  +\alpha\frac{225}{676},
 &w(\alpha)&=(1-\alpha)\frac{126}{1331}
                  +\alpha\frac{2025}{2197}.\label{eq:trialmoments}
\end{align}
Then $\E[u]=Ne$, while the mean, variance, and third central moment of $S$ are $Nm,Nv,Nw$, respectively. Indeed, a Bernoulli variable of success probability $p$ has variance $p(1-p)$ and third central moment $p(1-p)(1-2p)$. Multiplication by $3$ contributes factors $9$ and $27$ to these two moments, and centered moments of order two and three add for independent variables. These observations give every coefficient in \eqref{eq:trialmoments}.

Since $\D[u]=\mathbb E(S^{-1})$, Lemma~\ref{lem:harmonicmean} yields the upper bound
\begin{equation}\label{eq:trialpolynomial}
 \frac{C_{N,k}}{N^2}
 \le e(\alpha)m(\alpha)
       -\frac{e(\alpha)v(\alpha)^2}
                    {Nm(\alpha)v(\alpha)+w(\alpha)}.
\end{equation}
All four functions in \eqref{eq:trialmoments} are positive on $[0,1]$.

\emph{Step 2: the comparison for every $N\ge10$.}
For a real parameter $n\ge9$, set
\[
 G_n(\alpha)=\left(\alpha+\frac12-\frac1n\right)^2-em
                         +\frac{ev^2}{nmv+w},\qquad
 P_n(\alpha)=(nmv+w)G_n(\alpha),
\]
where $e,m,v,w$ are evaluated at $\alpha$. The function $P_n$ is a polynomial of degree at most four in $\alpha$. Positivity of $G_N(k/N)$ proves $C_{N,k}<A_{N,k}$ by \eqref{eq:trialpolynomial}.

Expansion of the affine functions \eqref{eq:trialmoments} gives the following coefficientwise bounds. If $P(t)=\sum_{j=0}^4c_jt^j$, the entries in its row are lower bounds for $c_0,\ldots,c_4$.
\[
\begin{array}{c|ccccc}
 P(t)&c_0&c_1&c_2&c_3&c_4\\ \hline
 P_{10}(t)&9/500&-13/1000&-27/100&7/12&0\\
 P_9(4/9+t)&1/250&9/50&1&7/5&1/2
\end{array}
\]
Appendix~\ref{app:trialcoefficients} lists the exact coefficients; comparison with the displayed fractions gives the stated bounds. The first row implies, for $\alpha\ge0$,
\[
 P_{10}(\alpha)\ge
 \frac9{500}-\frac{13}{1000}\alpha
 -\frac{27}{100}\alpha^2+\frac7{12}\alpha^3.
\]
The elementary square estimate
\[
 \frac{13}{1000}\alpha\le\frac1{50}\alpha^2+\frac{169}{80000}
\]
and the minimum of a cubic give
\[
 \frac7{12}\alpha^3-\frac{29}{100}\alpha^2
 \ge-\frac{4(29/100)^3}{27(7/12)^2}>-\frac{11}{1000}.
\]
Consequently,
\begin{equation}\label{eq:positivePten}
 P_{10}(\alpha)>\frac{391}{80000}>0
 \qquad(\alpha\ge0).
\end{equation}

It remains to show that the comparison persists as $n$ increases. The affine functions satisfy
\[
 v\le\frac9{11}m,\qquad
 \frac{ev}{m}\le\frac9{11}e
 =\frac{90}{121}+\frac{909}{1573}\alpha
 <\frac45+2\alpha\qquad(0\le\alpha\le1).
\]
For the first inequality, the two endpoint comparisons are $18/121=(9/11)(2/11)$ and $225/676<(9/11)(29/26)$, and both sides are affine. The last inequality follows by comparing the constant and linear terms. Since $w>0$, differentiation now gives
\begin{align*}
 \partial_nG_n(\alpha)
 &=\frac{2(\alpha+1/2-1/n)}{n^2}
      -\frac{emv^3}{(nmv+w)^2}\\
 &\ge\frac{2\alpha+1-2/n-ev/m}{n^2}>0
 \qquad(n\ge10,\ 0\le\alpha\le1).
\end{align*}
Thus $G_n\ge G_{10}>0$ on this range, proving the assertion for $N\ge10$ and $k\ge1$.

\emph{Step 3: dimension nine.}
If $k\ge4$, then $\alpha=k/9\ge4/9$. The second row of the coefficient comparison gives $P_9(\alpha)>0$, and the same product trial proves strictness. If $k=1$, Proposition~\ref{prop:strictcriterion} gives $C_{9,1}<18<81/4=A_{9,1}$.

For $k=2$, use instead the profiles
\[
 \widehat f=\frac{6\delta_1+\delta_2}{\sqrt{37}},\qquad
 \widehat g=\frac{\delta_{-1}+4\delta_0+\delta_1}{\sqrt{18}}.
\]
For their product with two half-line and seven line factors, the energy, mean, variance, and third central moment are
\[
 E=\frac{3706}{333},\quad M=\frac{979}{333},\quad
 V=\frac{648}{1369}+\frac{56}{81},\quad
 T=\frac{68040}{50653}+\frac{392}{729}.
\]
Substitution in \eqref{eq:harmonicmean} gives
\[
 C_{9,2}\le E\left(M-\frac{V^2}{MV+T}\right)
 =\frac{13519400682934}{452558671317}
 <30<\frac{121}{4}=A_{9,2}.
\]
The inequality with $30$ is equivalent to $57359456576>0$ after clearing the positive denominator.

For $k=3$, return to the original profiles $f,g$. Their product has energy $E=1473/143$. The identity $1/s=\int_0^1t^{s-1}\,dt$, applied to the law of $S=3+3I+J$, gives
\[
 D:=\D[u]=\frac1{26^3 11^6}
       \int_0^1t^2(25+t^3)^3(9+2t)^6\,dt.
\]
Retain only the constant and linear terms in $t^3$ from $(25+t^3)^3$; all omitted terms are positive. Expanding the remaining sixth power and integrating monomials yields
\begin{align*}
 D-\frac E{A_{9,3}}
 &\ge\frac1{26^3 11^6}\sum_{j=0}^6\binom6j9^{6-j}2^j
       \left(\frac{15625}{j+3}+\frac{1875}{j+6}\right)
       -\frac{5892}{24167}\\
 &=\frac{1212444977}{21577910602248}>0.
\end{align*}
Hence $C_{9,3}\le E/D<A_{9,3}$. This completes the mixed-sector proof.

\emph{Step 4: the punctured endpoint and optimality.}
The analytic punctured-lattice comparison in \cite[Proposition~2.1]{LizamaTransition} gives $C_{N,0}<A_{N,0}$ for every $N\ge9$, and $C_{8,0}=A_{8,0}$. The first assertion completes the strict comparison; the second prevents replacing $9$ by any smaller uniform threshold.
\end{proof}

The two elementary punctured cases in Table~\ref{tab:endpoints} follow from the half-line inequality and a cutoff argument, respectively. In dimension one, the classical discrete half-line Hardy inequality gives $C_{1,1}\ge1/4$; see \cite[Theorem~1.1]{LizamaPositive}. The dilation argument following Proposition~\ref{prop:continuum} gives the reverse inequality. The punctured line splits into two Dirichlet half-lines, with both forms splitting as sums, so $C_{1,0}=C_{1,1}=1/4$. In dimension two, choose a Lipschitz function $\chi$ that equals one on $[0,1]$ and zero on $[2,\infty)$, and put $u_R(x)=\chi(|x|_\infty/R)$ for $x\ne0$, with zero value at the origin. The four edges incident to the origin have bounded total energy. On all other edges the difference is $O(R^{-1})$, and only $O(R^2)$ edges contribute, so $\E[u_R]=O(1)$. The square shell $|x|_\infty=r$ has $8r$ points, each with $|x|^2\le2r^2$. For integer $R\ge1$, it follows that
\[
 \D[u_R]\ge4\sum_{r=1}^R\frac1r\longrightarrow\infty.
\]
Thus $C_{2,0}=0$. The positive-lattice equality $C_{2,2}=4$ and strict comparison for $N\ge3$ are proved in \cite[Theorem~1.1]{LizamaPositive}. 

\subsection{Critical inverse-square coupling}\label{sec:spectral}
The sharp Hardy coefficient gives the binding threshold for an attractive inverse-square potential. Let $J>0$, let $\mathsf{Adj}_{N,k}$ be the adjacency operator on $\Omega_{N,k}$, and let $W$ be multiplication by $|x|^{-2}$. Set
\begin{equation}\label{eq:Hg}
 \mathsf H_g=-J\mathsf{Adj}_{N,k}-gW,\qquad g\ge0.
\end{equation}
The boundary convention is the one in \eqref{eq:laplacian}, including the removed vertex when $k=0$.

\begin{corollary}\label{cor:spectral}
For $0\le k\le N$, the essential spectrum of $\mathsf H_g$ is $[-2NJ,2NJ]$. It has an eigenvalue below $-2NJ$ exactly when
\begin{equation}\label{eq:critical}
 g>JC_{N,k}.
\end{equation}
Thus the critical coupling $g_c(N,k)=JC_{N,k}$ satisfies the uniform law of Theorem~\ref{thm:main}, multiplied by $J$. In particular,
\[
 g_c(N,N)=J\kappa N^2+O(JN),\qquad
 g_c(N,0)=2JN+O(J\sqrt N).
\]
\end{corollary}
\begin{proof}
For $k\ge1$, Fourier and sine transforms identify $L_{N,k}$ with multiplication by $\sum_{i=1}^N(2-2\cos\theta_i)$, whose essential range is $[0,4N]$. For $k=0$, adjoining the removed vertex changes the full lattice operator by finite rank, giving the same essential spectrum. The operator $W$ is compact: its restriction to $|x|\le R$ has finite rank, and the remaining part has norm at most $R^{-2}$. Hence $\mathsf H_g=J(L_{N,k}-2NI)-gW$ has the stated essential spectrum.

If $g\le JC_{N,k}$, the Hardy inequality gives $\mathsf H_g+2NJI\ge0$. If $g>JC_{N,k}$, the defining infimum supplies a finitely supported $u$ with $J\E[u]-g\D[u]<0$. The variational principle then gives an eigenvalue below the essential band. The estimates in Propositions~\ref{prop:localization} and~\ref{prop:punctured} give the endpoint laws.
\end{proof}

At the critical coupling the shifted form is nonnegative. The classification of states at the band edge belongs to the criticality theory of graph Schr\"odinger operators; see \cite{KPPcriticality}.

\appendix
\section{An analytic comparison for the half-line oscillator}\label{app:comparison}

We prove Lemma~\ref{lem:comparison} by separating the first three coordinates from the tail. The only polynomial sign needed in the argument follows from concavity on a single interval.

\begin{proof}[Proof of Lemma~\ref{lem:comparison}]
\emph{Step 1: separate a finite block from the tail.}
Let $c=67/25$ and $1/3\le s\le1$. Begin with a finitely supported sequence $f$. At the edge joining coordinates $3$ and $4$, the inequality $2ab\le(2/3)a^2+(3/2)b^2$ gives
\[
 -2\Rea(\overline f_3f_4)\ge-\frac23|f_3|^2-\frac32|f_4|^2.
\]
The first three coordinates of $L_++s^2X^2-csI$ are then bounded below by the matrix
\begin{equation}\label{eq:Bmatrix}
 B(s)=\begin{pmatrix}
 2+s^2-cs&-1&0\\
 -1&2+4s^2-cs&-1\\
 0&-1&\frac43+9s^2-cs
 \end{pmatrix}.
\end{equation}
The Dirichlet kinetic form on the tail is $|f_4|^2+\sum_{n\ge4}|f_{n+1}-f_n|^2$. After subtracting $3|f_4|^2/2$, the tail is bounded below by
\begin{equation}\label{eq:tailpositive}
 \left(16s^2-cs-\frac12\right)I\ge\frac{173}{450}I.
\end{equation}
Here $16s^2-cs$ is increasing for $s\ge1/3$.

\emph{Step 2: reduce positivity of the block to its determinant.}
The lower right two-by-two block of $B(s)$ is positive definite, since
\[
 2+4s^2-cs\ge2-\frac{c^2}{16}=\frac{15511}{10000}>1,\qquad
 \frac43+9s^2-cs\ge\frac{36}{25}>1.
\]
The second bound follows because its left-hand side is increasing on $[1/3,1]$. By a Schur complement, it remains to prove $P(s)=\det B(s)>0$.

\emph{Step 3: prove the determinant sign on the whole interval.}
Write $a_1=2+s^2-cs$, $a_2=2+4s^2-cs$, and $a_3=4/3+9s^2-cs$. Expansion along the first row gives
\[
 P(s)=a_1(a_2a_3-1)-a_3.
\]
Subtracting the affine function with the same values at $s=1/3$ and $s=1$ leaves a polynomial vanishing at both endpoints. Factoring those two roots gives the exact identity
\begin{equation}\label{eq:Pidentity}
 P(s)=\frac{1-s}{2/3}P(1/3)+\frac{s-1/3}{2/3}P(1)
       +(s-1/3)(1-s)Q(s),
\end{equation}
where
\begin{gather}
 P(1/3)=\frac{2918}{46875},\qquad P(1)=\frac{2462}{15625},\label{eq:Pends}\\
 Q(s)=-\frac{562500s^4-1301875s^3+1137400s^2-555688s+93066}{15625}.
 \label{eq:Qpoly}
\end{gather}
Two differentiations and completion of the square give
\begin{equation}\label{eq:Qconcave}
 Q''(s)=-432\left(s-\frac{2083}{3600}\right)^2-\frac{28727}{30000}<0.
\end{equation}
Moreover,
\[
 Q(1/3)=\frac{38117}{84375}>0,\qquad Q(1)=\frac{64597}{15625}>0.
\]
A concave function lies above its chord, so $Q(s)>0$ throughout $[1/3,1]$. In \eqref{eq:Pidentity}, the first two terms form a strictly positive convex combination and the last term is nonnegative. Hence $P(s)>0$ and $B(s)$ is positive definite.

\emph{Step 4: recover the infinite-dimensional operator inequality.}
For each fixed $s$ in the interval, let $\mu(s)>0$ be the smallest eigenvalue of $B(s)$. The block and tail estimates give
\[
 \ip f{(L_++s^2X^2-csI)f}
 \ge\mu(s)\sum_{j=1}^3|f_j|^2
      +\frac{173}{450}\sum_{j\ge4}|f_j|^2
 \ge\min\{\mu(s),173/450\}\|f\|_2^2.
\]
Finite truncations are a core for the moment form domain, so this inequality extends to that domain. In particular the operator has a strictly positive lower bound, as asserted in \eqref{eq:comparison}. 
\end{proof}

\section{The small-parameter Mathieu coefficients}\label{app:coefficients}

For reference, we derive the two expansions used in Section~\ref{sec:smallalpha}. They agree with the standard tables \cite[Eqs.~28.6.1 and 28.6.5]{DLMF}, but the derivation below fixes their normalization directly from the matrices appearing in this paper.

Both matrices are of the form $D-qA_1$, with $D$ diagonal and $A_1$ bounded. At $q=0$ their ground eigenvalues are simple and isolated. Analytic perturbation theory therefore gives analytic eigenvalues and analytic eigenvectors near zero. We normalize each eigenvector by setting its first nonzero coordinate at $q=0$ equal to one for all small $q$. This is possible because that coordinate remains nonzero. Alternating the signs of successive coordinates changes $q$ to $-q$, so the eigenvalues are even functions of $q$.

We also record why distant coordinates have the orders used below. Expand the eigenvalue equation in powers of $q$. At order zero only the initial coordinate is nonzero. At each subsequent order, the perturbation couples only neighboring coordinates. Induction in the coefficient equations shows that a coordinate at distance $j$ from the initial one is $O(q^j)$. The diagonal gaps at $q=0$ are nonzero at every other coordinate, so each induction step can be solved. Analyticity makes these estimates genuine local remainder estimates.

\subsection{The even periodic eigenvalue}
In the even Fourier basis, indexed by $0,1,2,\ldots$, alternating signs turns the periodic Mathieu operator into a matrix with diagonal
\[
 0,4,16,36,\ldots
\]
and off-diagonal entries $-\sqrt2q,-q,-q,\ldots$. The factor $\sqrt2$ comes from the normalized constant function and normalized first cosine. The periodic ground eigenfunction is even near $q=0$, so its eigenvalue is $a=a_0(q)$.

Normalize its first coordinate to $v_0=1$. The first three coordinate equations are
\[
 a=-\sqrt2qv_1,\qquad
 (4-a)v_1=\sqrt2q+qv_2,\qquad
 (16-a)v_2=qv_1+qv_3.
\]
By analyticity, parity, and the distance estimate above, write
\[
 v_1=Aq+Bq^3+O(q^5),\qquad
 v_2=Cq^2+O(q^4),\qquad v_3=O(q^3),
\]
and $a=a_2q^2+a_4q^4+O(q^6)$. Comparing coefficients gives, successively,
\[
 4A=\sqrt2,\qquad a_2=-\sqrt2A=-\frac12,\qquad
 16C=A,\qquad 4B-a_2A=C.
\]
Thus
\[
 A=\frac{\sqrt2}{4},\qquad C=\frac{\sqrt2}{64},
 \qquad B=-\frac{7\sqrt2}{256},\qquad
 a_4=-\sqrt2B=\frac7{128}.
\]
This proves the first expansion in \eqref{eq:mathieusmall}.

\subsection{The Dirichlet eigenvalue}
For the matrix $T_q$ in \eqref{eq:Tq}, normalize the ground vector by $f_1=1$ and write $b=b_2(q)$. The first equations are
\[
 b=4-qf_2,\qquad
 (16-b)f_2=q+qf_3,\qquad
 (36-b)f_3=qf_2+qf_4.
\]
Use
\[
 f_2=Aq+Bq^3+O(q^5),\qquad f_3=Cq^2+O(q^4),\qquad
 f_4=O(q^3),
\]
and $b=4+\beta_2q^2+\beta_4q^4+O(q^6)$. The coefficient equations are
\[
 12A=1,\qquad \beta_2=-A=-\frac1{12},\qquad
 32C=A,\qquad 12B-\beta_2A=C.
\]
They give
\[
 A=\frac1{12},\qquad C=\frac1{384},\qquad
 B=-\frac5{13824},\qquad \beta_4=-B=\frac5{13824}.
\]

\begin{samepage}
\section{Exact coefficients for the strict-comparison trials}\label{app:trialcoefficients}

For the polynomials $P_n$ defined in the proof of Theorem~\ref{thm:comparison}, write $P_{10}(t)=\sum_{j=0}^4a_jt^j$ and $P_9(4/9+t)=\sum_{j=0}^4b_jt^j$. Their exact coefficients are
\begin{center}
\begingroup\renewcommand{\arraystretch}{2}
\begin{tabular}{@{}c@{\qquad}c@{\qquad}c@{}}
\toprule
$j$ & $a_j$ & $b_j$\\\midrule
$0$ & $\dfrac{73224}{4026275}$ & $\dfrac{15760971995}{3229043882922}$\\
$1$ & $-\dfrac{112718808}{8845726175}$ & $\dfrac{1567277401}{8491897128}$\\
$2$ & $-\dfrac{49117754871}{183991104440}$ & $\dfrac{1940652737663}{1913507486176}$\\
$3$ & $\dfrac{562894048809}{956753743088}$ & $\dfrac{1399054573629}{956753743088}$\\
$4$ & $\dfrac{280028354445}{478376871544}$ & $\dfrac{504051038001}{956753743088}$\\\bottomrule
\end{tabular}\endgroup
\end{center}

These coefficients exceed the corresponding lower bounds in Section~\ref{subsec:strictcomparison}.

\end{samepage}

\section{Numerical illustration}\label{app:illustration}
For Figure~\ref{fig:scaling}(a), $80$-dimensional Jacobi truncations have diagonal $4j^2$ ($j\ge0$) and off-diagonal $-\sqrt2q,-q,\ldots$ for $a_0$, and diagonal $4j^2$ ($j\ge1$) and off-diagonal $-q$ for $b_2$. SciPy~1.17.0's \texttt{eigh\_tridiagonal} computes their lowest eigenvalues with LAPACK driver \texttt{stebz} and tolerance $10^{-13}$.

At $\alpha=j/200$, $1\le j\le200$, we scan \eqref{eq:mathieupsi} on a $241$-point logarithmic grid in $q\in[10^{-5},10^4]$. Each interior grid minimum is refined between its neighbors in $\log q$ by bounded Brent minimization (SciPy's \texttt{minimize\_scalar}, \texttt{method='bounded'}, \texttt{xatol=1e-12}); the least candidate is retained. We set $\Psi(0)=0$. Increasing the truncation size to $160$ at the selected parameters changes the objective values by less than $10^{-9}$.

The coordinates plotted in Figure~\ref{fig:scaling}(a) are embedded in the LaTeX source.

\section*{Conflict of interest}
The corresponding author states that there is no conflict of interest.

\section*{Data availability}
The numerical data plotted in Figure~1(a) are included in the LaTeX source file accompanying this manuscript. The computational procedure is described in Appendix~D.

\section*{Declaration of generative AI and AI-assisted technologies in the manuscript preparation process}
I used OpenAI's ChatGPT for mathematical exploration, checks of candidate arguments and calculations, literature searches, and drafting. I reviewed and edited the material and take responsibility for the mathematical claims, references, and final text.

\end{document}